
%
%
%
\def\unredoffs{} \def\redoffs{\voffset=-.31truein\hoffset=-.59truein}
\def\speclscape{\special{ps: landscape}}
%
%
%
%
\newbox\leftpage \newdimen\fullhsize \newdimen\hstitle \newdimen\hsbody
\tolerance=1000\hfuzz=2pt
\catcode`\@=11 
%
\ifx\answ\bigans\message{(This will come out unreduced.}
\magnification=1200\unredoffs\baselineskip=16pt plus 2pt minus 1pt
\hsbody=\hsize \hstitle=\hsize 
\else\message{(This will be reduced.} \let\l@r=L
\magnification=1000\baselineskip=16pt plus 2pt minus 1pt \vsize=7truein
\redoffs \hstitle=8truein\hsbody=4.75truein\fullhsize=10truein\hsize=\hsbody
\output={\ifnum\pageno=0 
  \shipout\vbox{\speclscape{\hsize\fullhsize\makeheadline}
    \hbox to \fullhsize{\hfill\pagebody\hfill}}\advancepageno
  \else
  \almostshipout{\leftline{\vbox{\pagebody\makefootline}}}\advancepageno
  \fi}
\def\almostshipout#1{\if L\l@r \count1=1 \message{[\the\count0.\the\count1]}
      \global\setbox\leftpage=#1 \global\let\l@r=R
 \else \count1=2
  \shipout\vbox{\speclscape{\hsize\fullhsize\makeheadline}
      \hbox to\fullhsize{\box\leftpage\hfil#1}}  \global\let\l@r=L\fi}
\fi
%
\newcount\yearltd\yearltd=\year\advance\yearltd by -1900

%

%
%

\def\draftmode{\message{ DRAFTMODE }\def\draftdate{{\rm preliminary draft:
\number\month/\number\day/\number\yearltd\ \ \hourmin}}%
\headline={\hfil\draftdate}\writelabels\baselineskip=20pt plus 2pt minus 2pt
 {\count255=\time\divide\count255 by 60 \xdef\hourmin{\number\count255}
  \multiply\count255 by-60\advance\count255 by\time
  \xdef\hourmin{\hourmin:\ifnum\count255<10 0\fi\the\count255}}}
\def\nolabels{\def\wrlabeL##1{}\def\eqlabeL##1{}\def\reflabeL##1{}}
\def\writelabels{\def\wrlabeL##1{\leavevmode\vadjust{\rlap{\smash%
{\line{{\escapechar=` \hfill\rlap{\sevenrm\hskip.03in\string##1}}}}}}}%
\def\eqlabeL##1{{\escapechar-1\rlap{\sevenrm\hskip.05in\string##1}}}%
\def\reflabeL##1{\noexpand\llap{\noexpand\sevenrm\string\string\string##1}}}
\nolabels
%
\global\newcount\secno \global\secno=0
\global\newcount\meqno \global\meqno=1
\def\newsec#1{\global\advance\secno by1\message{(\the\secno. #1)}
\global\subsecno=0\eqnres@t\noindent{\bf\the\secno. #1}
\writetoca{{\secsym} {#1}}\par\nobreak\medskip\nobreak}
\def\eqnres@t{\xdef\secsym{\the\secno.}\global\meqno=1\bigbreak\bigskip}
\def\sequentialequations{\def\eqnres@t{\bigbreak}}\xdef\secsym{}
\global\newcount\subsecno \global\subsecno=0
\def\subsec#1{\global\advance\subsecno by1\message{(\secsym\the\subsecno. #1)}
\ifnum\lastpenalty>9000\else\bigbreak\fi
\noindent{\it\secsym\the\subsecno. #1}\writetoca{\string\quad
{\secsym\the\subsecno.} {#1}}\par\nobreak\medskip\nobreak}
\def\appendix#1#2{\global\meqno=1\global\subsecno=0\xdef\secsym{\hbox{#1.}}
\bigbreak\bigskip\noindent{\bf Appendix #1. #2}\message{(#1. #2)}
\writetoca{Appendix {#1.} {#2}}\par\nobreak\medskip\nobreak}
%
%
\def\eqnn#1{\xdef #1{(\secsym\the\meqno)}\writedef{#1\leftbracket#1}%
\global\advance\meqno by1\wrlabeL#1}
\def\eqna#1{\xdef #1##1{\hbox{$(\secsym\the\meqno##1)$}}
\writedef{#1\numbersign1\leftbracket#1{\numbersign1}}%
\global\advance\meqno by1\wrlabeL{#1$\{\}$}}
\def\eqn#1#2{\xdef #1{(\secsym\the\meqno)}\writedef{#1\leftbracket#1}%
\global\advance\meqno by1$$#2\eqno#1\eqlabeL#1$$}
%
\newskip\footskip\footskip14pt plus 1pt minus 1pt 
\def\footnotefont{\ninepoint}\def\f@t#1{\footnotefont #1\@foot}
\def\f@@t{\baselineskip\footskip\bgroup\footnotefont\aftergroup\@foot\let\next}
\setbox\strutbox=\hbox{\vrule height9.5pt depth4.5pt width0pt}
\global\newcount\ftno \global\ftno=0
\def\foot{\global\advance\ftno by1\footnote{$^{\the\ftno}$}}
%
\newwrite\ftfile
\def\footend{\def\foot{\global\advance\ftno by1\chardef\wfile=\ftfile
$^{\the\ftno}$\ifnum\ftno=1\immediate\openout\ftfile=foots.tmp\fi%
\immediate\write\ftfile{\noexpand\smallskip%
\noexpand\item{f\the\ftno:\ }\pctsign}\findarg}%
\def\footatend{\vfill\eject\immediate\closeout\ftfile{\parindent=20pt
\centerline{\bf Footnotes}\nobreak\bigskip\input foots.tmp }}}
\def\footatend{}
%
%
\global\newcount\refno \global\refno=1
\newwrite\rfile
%
\def\ref{\nref}
\def\nref#1{\xdef#1{[\the\refno]}\writedef{#1\leftbracket#1}%
\ifnum\refno=1\immediate\openout\rfile=refs.tmp\fi
\global\advance\refno by1\chardef\wfile=\rfile\immediate
\write\rfile{\noexpand\item{#1\ }\reflabeL{#1\hskip.31in}\pctsign}\findarg}
\def\findarg#1#{\begingroup\obeylines\newlinechar=`\^^M\pass@rg}
{\obeylines\gdef\pass@rg#1{\writ@line\relax #1^^M\hbox{}^^M}%
\gdef\writ@line#1^^M{\expandafter\toks0\expandafter{\striprel@x #1}%
\edef\next{\the\toks0}\ifx\next\em@rk\let\next=\endgroup\else\ifx\next\empty%
\else\immediate\write\wfile{\the\toks0}\fi\let\next=\writ@line\fi\next\relax}}
\def\striprel@x#1{} \def\em@rk{\hbox{}}
\def\lref{\begingroup\obeylines\lr@f}
\def\lr@f#1#2{\gdef#1{\ref#1{#2}}\endgroup\unskip}

\def\addref#1{\immediate\write\rfile{\noexpand\item{}#1}} 
\def\startrefs#1{\immediate\openout\rfile=refs.tmp\refno=#1}
\def\refs#1{\count255=1[\r@fs #1{\hbox{}}]}
\def\r@fs#1{\ifx\und@fined#1\message{reflabel \string#1 is undefined.}%
\nref#1{need to supply reference \string#1.}\fi%
\vphantom{\hphantom{#1}}\edef\next{#1}\ifx\next\em@rk\def\next{}%
\else\ifx\next#1\ifodd\count255\relax\xref#1\count255=0\fi%
\else#1\count255=1\fi\let\next=\r@fs\fi\next}
%

%
\newwrite\ffile\global\newcount\figno \global\figno=1
\def\fig{fig.~\the\figno\nfig}
\def\nfig#1{\xdef#1{fig.~\the\figno}%
\writedef{#1\leftbracket fig.\noexpand~\the\figno}%
\ifnum\figno=1\immediate\openout\ffile=figs.tmp\fi\chardef\wfile=\ffile%
\immediate\write\ffile{\noexpand\medskip\noexpand\item{Fig.\ \the\figno. }
\reflabeL{#1\hskip.55in}\pctsign}\global\advance\figno by1\findarg}
\def\xfig{\expandafter\xf@g}\def\xf@g fig.\penalty\@M\ {}
\def\figs#1{figs.~\f@gs #1{\hbox{}}}
\def\f@gs#1{\edef\next{#1}\ifx\next\em@rk\def\next{}\else
\ifx\next#1\xfig #1\else#1\fi\let\next=\f@gs\fi\next}
\newwrite\lfile
{\escapechar-1\xdef\pctsign{\string\%}\xdef\leftbracket{\string\{}
\xdef\rightbracket{\string\}}\xdef\numbersign{\string\#}}

\def\writestop{\def\writestoppt{\immediate\write\lfile{\string\pageno%
\the\pageno\string\startrefs\leftbracket\the\refno\rightbracket%
\string\def\string\secsym\leftbracket\secsym\rightbracket%
\string\secno\the\secno\string\meqno\the\meqno}\immediate\closeout\lfile}}
\def\writestoppt{}\def\writedef#1{}
\def\seclab#1{\xdef #1{\the\secno}\writedef{#1\leftbracket#1}\wrlabeL{#1=#1}}
\def\subseclab#1{\xdef #1{\secsym\the\subsecno}%
\writedef{#1\leftbracket#1}\wrlabeL{#1=#1}}
\newwrite\tfile \def\writetoca#1{}
\def\leaderfill{\leaders\hbox to 1em{\hss.\hss}\hfill}
\def\writetoc{\immediate\openout\tfile=toc.tmp
   \def\writetoca##1{{\edef\next{\write\tfile{\noindent ##1
   \string\leaderfill {\noexpand\number\pageno} \par}}\next}}}
%
%
%

\catcode`\@=12 
%
\edef\tfontsize{\ifx\answ\bigans scaled\magstep3\else scaled\magstep4\fi}
\font\titlerm=cmr10 \tfontsize \font\titlerms=cmr7 \tfontsize
\font\titlermss=cmr5 \tfontsize \font\titlei=cmmi10 \tfontsize
\font\titleis=cmmi7 \tfontsize \font\titleiss=cmmi5 \tfontsize
\font\titlesy=cmsy10 \tfontsize \font\titlesys=cmsy7 \tfontsize
\font\titlesyss=cmsy5 \tfontsize \font\titleit=cmti10 \tfontsize
\skewchar\titlei='177 \skewchar\titleis='177 \skewchar\titleiss='177
\skewchar\titlesy='60 \skewchar\titlesys='60 \skewchar\titlesyss='60
\def\titlefont{\def\rm{\fam0\titlerm}
\textfont0=\titlerm \scriptfont0=\titlerms \scriptscriptfont0=\titlermss
\textfont1=\titlei \scriptfont1=\titleis \scriptscriptfont1=\titleiss
\textfont2=\titlesy \scriptfont2=\titlesys \scriptscriptfont2=\titlesyss
\textfont\itfam=\titleit \def\it{\fam\itfam\titleit}\rm}
 \ifx\answ\bigans\else scaled\magstep1\fi
\ifx\answ\bigans\else

 \font\absi=cmmi10 scaled\magstep1
\font\absis=cmmi7 scaled\magstep1 \font\absiss=cmmi5 scaled\magstep1
\font\abssy=cmsy10 scaled\magstep1 \font\abssys=cmsy7 scaled\magstep1
\font\abssyss=cmsy5 scaled\magstep1 
\skewchar\absi='177 \skewchar\absis='177 \skewchar\absiss='177
\skewchar\abssy='60 \skewchar\abssys='60 \skewchar\abssyss='60
\fi
\font\ninerm=cmr9 \font\sixrm=cmr6 \font\ninei=cmmi9 \font\sixi=cmmi6
\font\ninesy=cmsy9 \font\sixsy=cmsy6 \font\ninebf=cmbx9
\font\nineit=cmti9 \font\ninesl=cmsl9 \skewchar\ninei='177
\skewchar\sixi='177 \skewchar\ninesy='60 \skewchar\sixsy='60
\def\ninepoint{\def\rm{\fam0\ninerm}
\textfont0=\ninerm \scriptfont0=\sixrm \scriptscriptfont0=\fiverm
\textfont1=\ninei \scriptfont1=\sixi \scriptscriptfont1=\fivei
\textfont2=\ninesy \scriptfont2=\sixsy \scriptscriptfont2=\fivesy
\textfont\itfam=\ninei \def\it{\fam\itfam\nineit}\def\sl{\fam\slfam\ninesl}%
\textfont\bffam=\ninebf \def\bf{\fam\bffam\ninebf}\rm}
%
%

\hyphenation{anom-aly anom-alies coun-ter-term coun-ter-terms}
\def\inv{^{\raise.15ex\hbox{${\scriptscriptstyle -}$}\kern-.05em 1}}

\def\Dsl{\,\raise.15ex\hbox{/}\mkern-13.5mu D} 
\def\dsl{\raise.15ex\hbox{/}\kern-.57em\partial}

\def\lspace{\ifx\answ\bigans{}\else\qquad\fi}
\def\lbspace{\ifx\answ\bigans{}\else\hskip-.2in\fi} 
\def\boxeqn#1{\vcenter{\vbox{\hrule\hbox{\vrule\kern3pt\vbox{\kern3pt
    \hbox{${\displaystyle #1}$}\kern3pt}\kern3pt\vrule}\hrule}}}
\def\mbox#1#2{\vcenter{\hrule \hbox{\vrule height#2in
        \kern#1in \vrule} \hrule}}  
%

\def\darr#1{\raise1.5ex\hbox{$\leftrightarrow$}\mkern-16.5mu #1}

\def\roughly#1{\raise.3ex\hbox{$#1$\kern-.75em\lower1ex\hbox{$\sim$}}}

%
%


\def\frac#1#2{{#1\over#2}}

\def\journal#1&#2(#3){\unskip, #1~\bf #2 \rm(19#3) }
\def\andjournal#1&#2(#3){\sl #1~\bf #2 \rm (19#3) }

\def\bra#1{\left\langle #1\right|}
\def\ket#1{\left| #1\right\rangle}

\catcode`\@=11\def\slash#1{\mathord{\mathpalette\c@ncel{#1}}}
\overfullrule=0pt
\def\steepslash{\c@ncel}
\def\frac#1#2{{#1\over #2}}

\def\:{\!:\!}
\def\inbar{\,\vrule height1.5ex width.4pt depth0pt}
\def\IQ{\relax\,\hbox{$\inbar\kern-.3em{\rm Q}$}}
\def\IB{\relax{\rm I\kern-.18em B}}
\def\IC{\relax\hbox{$\inbar\kern-.3em{\rm C}$}}
\def\IP{\relax{\rm I\kern-.18em P}}
\def\IR{\relax{\rm I\kern-.18em R}}
\def\ZZ{\relax\ifmmode\mathchoice
{\hbox{Z\kern-.4em Z}}{\hbox{Z\kern-.4em Z}}
{\lower.9pt\hbox{Z\kern-.4em Z}}
{\lower1.2pt\hbox{Z\kern-.4em Z}}\else{Z\kern-.4em Z}\fi}

\catcode`\@=12

\def\npb#1(#2)#3{{ Nucl. Phys. }{B#1} (#2) #3}
\def\plb#1(#2)#3{{ Phys. Lett. }{#1B} (#2) #3}
\def\pla#1(#2)#3{{ Phys. Lett. }{#1A} (#2) #3}
\def\prl#1(#2)#3{{ Phys. Rev. Lett. }{#1} (#2) #3}
\def\mpla#1(#2)#3{{ Mod. Phys. Lett. }{A#1} (#2) #3}
\def\ijmpa#1(#2)#3{{ Int. J. Mod. Phys. }{A#1} (#2) #3}
\def\cmp#1(#2)#3{{ Comm. Math. Phys. }{#1} (#2) #3}
\def\cqg#1(#2)#3{{ Class. Quantum Grav. }{#1} (#2) #3}
\def\jmp#1(#2)#3{{ J. Math. Phys. }{#1} (#2) #3}
\def\anp#1(#2)#3{{ Ann. Phys. }{#1} (#2) #3}
\def\prd#1(#2)#3{{ Phys. Rev. } {D{#1}} (#2) #3}
\def\ptp#1(#2)#3{{ Progr. Theor. Phys. }{#1} (#2) #3}
\def\aom#1(#2)#3{{ Ann. Math. }{#1} (#2) #3}

\def\bs{\bigskip}

\def\br{\buildrel}
\def\bra{\langle}
\def\ket{\rangle}

\def\C{{\bf C}}

\def\Z{{\bf Z}}
\def\cA{{\cal A}}
\def\cB{{\cal B}}
\def\cC{{\cal C}}
\def\cD{{\cal D}}
\def\cE{{\cal E}}
\def\cF{{\cal F}}

\def\cH{{\cal H}}
\def\cI{{\cal I}}

\def\cM{{\cal M}}

\def\cO{{\cal O}}

\def\cS{{\cal S}}

\def\cV{{\cal V}}
\def\cW{{\cal W}}

\input amssym

\def\gg{{\goth g}}

\def\gU{{\goth U}}

\def\da#1{{\p \over \p s_{#1}}}

\def\cicy#1(#2|#3)#4{\left(\matrix{#2}\right|\!\!
                     \left|\matrix{#3}\right)^{{#4}}_{#1}}

\def\ra{\rightarrow}

\def\da{\downarrow}

\def\bs{\bigskip}

\def\Box{{\,\lower0.9pt\vbox{\hrule
\hbox{\vrule height 0.2 cm \hskip 0.2 cm
\vrule height 0.2 cm}\hrule}\,}}

\global\newcount\thmno \global\thmno=0
\def\definition#1{\global\advance\thmno by1
\bigskip\noindent{\bf Definition \secsym\the\thmno. }{\it #1}
\par\nobreak\medskip\nobreak}
\def\question#1{\global\advance\thmno by1
\bigskip\noindent{\bf Question \secsym\the\thmno. }{\it #1}
\par\nobreak\medskip\nobreak}
\def\theorem#1{\global\advance\thmno by1
\bigskip\noindent{\bf Theorem \secsym\the\thmno. }{\it #1}
\par\nobreak\medskip\nobreak}
\def\proposition#1{\global\advance\thmno by1
\bigskip\noindent{\bf Proposition \secsym\the\thmno. }{\it #1}
\par\nobreak\medskip\nobreak}
\def\corollary#1{\global\advance\thmno by1
\bigskip\noindent{\bf Corollary \secsym\the\thmno. }{\it #1}
\par\nobreak\medskip\nobreak}
\def\lemma#1{\global\advance\thmno by1
\bigskip\noindent{\bf Lemma \secsym\the\thmno. }{\it #1}
\par\nobreak\medskip\nobreak}
\def\conjecture#1{\global\advance\thmno by1
\bigskip\noindent{\bf Conjecture \secsym\the\thmno. }{\it #1}
\par\nobreak\medskip\nobreak}
\def\exercise#1{\global\advance\thmno by1
\bigskip\noindent{\bf Exercise \secsym\the\thmno. }{\it #1}
\par\nobreak\medskip\nobreak}
\def\remark#1{\global\advance\thmno by1
\bigskip\noindent{\bf Remark \secsym\the\thmno. }{\it #1}
\par\nobreak\medskip\nobreak}
\def\problem#1{\global\advance\thmno by1
\bigskip\noindent{\bf Problem \secsym\the\thmno. }{\it #1}
\par\nobreak\medskip\nobreak}
\def\others#1#2{\global\advance\thmno by1
\bigskip\noindent{\bf #1 \secsym\the\thmno. }{\it #2}
\par\nobreak\medskip\nobreak}
\def\proof{\noindent Proof: }

\def\thmlab#1{\xdef #1{\secsym\the\thmno}\writedef{#1\leftbracket#1}\wrlabeL{#1=#1}}
%
%
\def\newsec#1{\global\advance\secno by1\message{(\the\secno. #1)}
\global\subsecno=0\thmno=0\eqnres@t\noindent{\bf\the\secno. #1}
\writetoca{{\secsym} {#1}}\par\nobreak\medskip\nobreak}
\def\eqnres@t{\xdef\secsym{\the\secno.}\global\meqno=1\bigbreak\bigskip}
\def\sequentialequations{\def\eqnres@t{\bigbreak}}\xdef\secsym{}
%

%
\newcount{\exnum}
\def\prob{\advance\exnum by 1
\bigskip\item{\the\exnum.}\ }
\newcount{\exnum}
\def\next{\advance\exnum by 1
\bigskip\noindent{\the\exnum.}\ }


\ref\BMP{P. Bouwknegt, J. McCarthy, K. Pilch, The $\cW_3$ algebra: Modules, Semi-infinite Cohomology and BV-algebras, Lect. Notes in Phys, New Series Monographs 42, Springer-Verlag, 1996.}
\ref\FF{ B. Feigin and E. Frenkel,   Semi-Infinite Weil complex and the Virasoro algebra,  Comm. Math. Phys. 137 (1991), 617-639.}
\ref\FI{ I. Frenkel, Two constructions of affine Lie algebras and boson-fermion correspondence in quantum field theory, J. Funct. Anal. 44 (1981) 259-327.}
\ref\FZ{ I.B. Frenkel, and Y.C. Zhu, Vertex operator algebras associated to representations of affine and Virasoro algebras, Duke Mathematical Journal, Vol. 66, No. 1, (1992), 123-168.}
\ref\FMS{ D. Friedan, E. Martinec, S. Shenker,   Conformal invariance, supersymmetry and string theory, Nucl. Phys. B271 (1986) 93-165.} 
\ref\GKO{ P. Goddard, A. Kent, and D. Olive, Virasoro algebras and coset space models,  Phys. Lett B 152 (1985) 88-93.}
\ref\GMS{V. Gorbounov, F. Malikov, V. Schectman, Gerbes of chiral differential operators, Math. Res. Lett. 7 (2000), no. 1, 55--66.}
\ref\HCI{Harish-Chandra, Differential operators on a semisimple Lie algebra, Am. J. Math, Vol. 79, No. 1 (1957) 87-120.}
\ref\HCII{Harish-Chandra, Invariant differential operators and distributions on a semisimple Lie algebra, Am. J. Math. Vol. 86, No. 3 (1964) 534-564.}
\ref\KP{V. Kac, D. Peterson, Infinite-dimensional Lie algebras, theta functions and modular forms, Adv. Math. 53 (1984) 125-264.}
\ref\KR{V. Kac, A. Radul, Representation theory of the vertex algebra $\cW_{1+\infty}$, Transf. Groups, Vol 1 (1996) 41-70.}
 \ref\Kn{F. Knop, A Harish-Chandra homomorphism for reductive group actions, Ann. Math. 140 (1994), 253-288.}
\ref\LL{B. Lian, A. Linshaw, Howe pairs in the theory of vertex algebras, J. Algebra 317, 111-152 (2007).}
\ref\LLI{B. Lian, and A. Linshaw, Chiral equivariant cohomology I, Adv. Math. 209, 99-161 (2007).}
\ref\MSV{F. Malikov, V. Schectman, and A. Vaintrob, Chiral de Rham complex, Commun. Math. Phys, 204, (1999) 439-473.}
\ref\Mu{I. Musson, M. van den Bergh, Invariants under Tori of Rings of Invariant Operators and Related Topics, Mem. Am. Math. Soc. No. 650 (1998).}
\ref\Sch{G. Schwarz, Finite-dimensional representations of invariant differential operators, J. Algebra 258 (2002) 160-204.}
\ref\WI{W. Wang, $\cW_{1+\infty}$ algebra, $\cW_3$ algebra, and Friedan-Martinec-Shenker bosonization, Commun. Math. Phys. 195 (1998), 95--111.}
\ref\WII{W. Wang, Classification of irreducible modules of $\cW_3$ with $c=-2$, Commun. Math. Phys. 195 (1998), 113--128.}
\ref\Za{A. Zamolodchikov, Infinite additional symmetries in two-dimensional conformal field theory, Theor. Math. Phys. 65 (1985) 1205-1213.}
\ref\Zh{Y. Zhu, Modular invariants of characters of vertex operators, J. Amer. Math. Soc. 9 (1996) 237-302.}

\centerline{\titlefont Invariant chiral differential operators}

\centerline{\titlefont and the $\cW_3$ algebra}

\bs
\centerline{Andrew R. Linshaw}
\bs

\baselineskip=13pt plus 2pt minus 2pt
ABSTRACT.  
Attached to a vector space $V$ is a vertex algebra $\cS(V)$ known as the $\beta\gamma$-system or algebra of chiral differential operators on $V$. It is analogous to the Weyl algebra $\cD(V)$, and is related to $\cD(V)$ via the Zhu functor. If $G$ is a connected Lie group with Lie algebra $\gg$, and $V$ is a linear $G$-representation, there is an action of the corresponding affine algebra on $\cS(V)$. The invariant space $\cS(V)^{\gg[t]}$ is a commutant subalgebra of $\cS(V)$, and plays the role of the classical invariant ring $\cD(V)^{G}$. When $G$ is an abelian Lie group acting diagonally on $V$, we find a finite set of generators for $\cS(V)^{\gg[t]}$, and show that $\cS(V)^{\gg[t]}$ is a simple vertex algebra and a member of a Howe pair. The Zamolodchikov $\cW_3$ algebra with $c=-2$ plays a fundamental role in the structure of $\cS(V)^{\gg[t]}$. 

\baselineskip=15pt plus 2pt minus 1pt
\parskip=\baselineskip

\centerline{\bf Contents}\nobreak\medskip{\baselineskip=12pt
 \parskip=0pt\catcode`\@=11 \input toc.tmp \catcode`\@=12 \bigbreak\bigskip}  

\newsec{Introduction}

Let $G$ be a connected, reductive Lie group acting algebraically on a smooth variety $X$. Throughout this paper, our base field will always be $\C$. The ring $\cD(X)^G$ of invariant differential operators on $X$ has been much studied in recent years. In the case where $X$ is the homogeneous space $G/K$, $\cD(X)^G$ was originally studied by Harish-Chandra in order to understand the various function spaces attached to $X$ \HCI\HCII. In general, $\cD(X)^G$ is not a homomorphic image of the universal enveloping algebra of a Lie algebra, but it is believed that $\cD(X)^G$ shares many properties of enveloping algebras. For example, the center of $\cD(X)^G$ is always a polynomial ring \Kn. In the case where $G$ is a torus, the structure and representation theory of the rings $\cD(X)^G$ were studied extensively in \Mu, but much less is known about $\cD(X)^G$ when $G$ is nonabelian. The first step in this direction was taken by Schwarz in \Sch, in which he considered the special but nontrivial case where $G = SL(3)$ and $X$ is the adjoint representation. In this case, he found generators for $\cD(X)^G$, showed that $\cD(X)^G$ is an FCR algebra, and classified its finite-dimensional modules.

\subsec{A vertex algebra analogue of $\cD(X)^G$}
In \MSV, Malikov-Schechtman-Vaintrob introduced a sheaf of vertex algebras on any smooth variety $X$ known as the chiral de Rham complex. For an affine open set $V\subset X$, the algebra of sections over $V$ is just a copy of the $bc\beta\gamma$-system $\cS(V)\otimes \cE(V)$, localized over the function ring $\cO(V)$. A natural question is whether there exists a subsheaf of \lq\lq chiral differential operators" on $X$, whose space of sections over $V$ is just the (localized) $\beta\gamma$-system $\cS(V)$. For general $X$, there is a cohomological obstruction to the existence of such a sheaf, but it does exist in certain special cases such as affine spaces and certain homogeneous spaces \MSV\GMS. 

In this paper, we focus on the case where $X$ is the affine space $V = \C^n$, and we take $\cS(V)$ to be our algebra of chiral differential operators on $V$. $\cS(V)$ is related to $\cD(V)$ via the {\it Zhu functor}, which attaches to every vertex algebra $\cV$ an associative algebra $A(\cV)$ known as the {\it Zhu algebra} of $\cV$, together with a surjective linear map $\pi_{Zh}:\cV\ra A(\cV)$. 

If $V$ carries a linear action of a group $G$ with Lie algebra $\gg$, the corresponding representation $\rho:\gg\ra End(V)$ induces a vertex algebra homomorphism \eqn\hatrho{\cO(\gg,B)\ra \cS(V).} Here $\cO(\gg,B)$ is the current algebra of $\gg$ associated to the bilinear form $B(\xi,\eta) = -Tr(\rho(\xi)\rho(\eta))$ on $\gg$. Letting $\Theta$ denote the image of $\cO(\gg,B)$ inside $\cS(V)$, the commutant $Com(\Theta,\cS(V))$, which we denote by $\cS(V)^{\Theta_+}$, is just the invariant space $\cS(V)^{\gg[t]}$. Accordingly, we call $\cS(V)^{\Theta_+}$ the algebra of {\it invariant chiral differential operators} on $V$. There is a commutative diagram
\eqn\commdiagva{\matrix{\cS(V)^{\Theta_+}&\br\iota\over\hookrightarrow&\cS(V)\cr \pi\da\hskip.4in& &\pi_{Zh}\da\hskip.4in\cr \cD(V)^{G}&\br\iota\over\hookrightarrow &\cD(V)}.}
Here the horizontal maps are inclusions, and the map $\pi$ on the left is the restriction of the Zhu map on $\cS(V)$ to the subalgebra $\cS(V)^{\Theta_+}$. In general, $\pi$ is not surjective, and $\cD(V)^{G}$ need not be the Zhu algebra of $\cS(V)^{\Theta_+}$. 

For a general vertex algebra $\cV$ and subalgebra $\cA$, the commutant $Com(\cA,\cV)$ was introduced by Frenkel-Zhu in \FZ, generalizing a previous construction in representation theory \KP~and conformal field theory \GKO~known as the coset construction. We regard $\cV$ as a module over $\cA$ via the left regular action, and we regard $Com(\cA,\cV)$, which we often denote by $\cV^{\cA_+}$, as the invariant subalgebra. Finding a set of generators for $\cV^{\cA_+}$, or even determining when it is finitely generated as a vertex algebra, is generally a non-trivial problem. It is also natural to study the double commutant $Com(\cV^{\cA_+},\cV)$, which always contains $\cA$. If $\cA = Com(\cV^{\cA_+},\cV)$, we say that $\cA$ and $\cV^{\cA_+}$ form a {\it Howe pair} inside $\cV$. Since $$Com(Com(\cV^{\cA_+},\cV),\cV) = \cV^{\cA_+},$$ a subalgebra $\cB$ is a member of a Howe pair if and only if $\cB = \cV^{\cA_+}$ for some $\cA$.

Here are some natural questions one can ask about $\cS(V)^{\Theta_+}$ and its relationship to $\cD(V)^{G}$.

\others{Question}{When is $\cS(V)^{\Theta_+}$ finitely generated as a vertex algebra? Can we find a set of generators?} 

\others{Question}{When do $\cS(V)^{\Theta_+}$ and $\Theta$ form a Howe pair inside $\cS(V)$? In the case where $G = SL(2)$ and $V$ is the adjoint module, this question was answered affirmatively in \LL.}

\others{Question}{What are the vertex algebra ideals in $\cS(V)^{\Theta_+}$, and when is $\cS(V)^{\Theta_+}$ a simple vertex algebra? }

\others{Question}{When is $\cS(V)^{\Theta_+}$ a conformal vertex algebra?}

\others{Question}{When is $\pi: \cS(V)^{\Theta_+}\ra \cD(V)^{G}$ surjective? More generally, describe $Im(\pi)$ and $Coker(\pi)$.}  

These questions are somewhat outside the realm of classical invariant theory because the Lie algebra $\gg[t]$ is both infinite-dimensional and non-reductive. Moreover, when $G$ is nonabelian, $\cS(V)$ need not decompose into a sum of irreducible $\cO(\gg,B)$-modules. The case where $G$ is simple and $V$ is the adjoint module is of particular interest to us, since in this case $\cS(V)^{\Theta_+}$ is a subalgebra of the complex $(\cW(\gg)_{bas},d)$ which computes the chiral equivariant cohomology of a point \LLI. 

In this paper, we focus on the case where $G$ is an abelian group acting faithfully and diagonalizably on $V$. This is much easier than the general case because $\cO(\gg,B)$ is then a tensor product of Heisenberg vertex algebras, which act completely reducibly on $\cS(V)$. For any such action, we find a finite set of generators for $\cS(V)^{\Theta_+}$, and show that $\cS(V)^{\Theta_+}$ is a simple vertex algebra. Moreover, $\cS(V)^{\Theta_+}$ and $\Theta$ always form a Howe pair inside $\cS(V)$. For generic actions, we show that $\cS(V)^{\Theta_+}$ admits a $k$-parameter family of conformal structures where $k = dim~V-dim~\gg$, and we find a finite set of generators for $Im(\pi)$. Finally, we show that $Coker(\pi)$ is always a finitely generated module over $Im(\pi)$ with generators corresponding to central elements of $\cD(V)^{G}$. The Zamolodchikov $\cW_3$ algebra of central charge $c=-2$ plays an important role in the structure of $\cS(V)^{\Theta_+}$. Our description relies on the fundamental papers \WI~\WII~of W. Wang, in which he classified the irreducible modules of $\cW_{3,-2}$.

In the case where $G$ is nonabelian, very little is known about the structure of $\cS(V)^{\Theta_+}$, and the representation-theoretic techniques used in the abelian case cannot be expected to work. In a separate paper, we will use tools from commutative algebra to describe $\cS(V)^{\Theta_+}$ in the special cases where $G$ is one of the classical Lie groups $SL(n)$, $SO(n)$, or $Sp(2n)$, and $V$ is a direct sum of copies of the standard representation.

One hopes that the vertex algebra point of view can also shed some light on the classical algebras $\cD(V)^{G}$. For example, the vertex algebra products on $\cS(V)$ induce a family of bilinear operations $*_k$, $k\geq -1$ on $\cD(V)^{G}$, which coincide with classical operations known as transvectants. $\cD(V)^{G}$ is generally not simple as an associative algebra, but in the case where $G$ is an abelian group acting diagonalizably on $V$, $\cD(V)^{G}$ is always simple as a $*$-algebra in the obvious sense.

\subsec{Acknowledgements} I thank B. Lian for helpful conversations and for suggesting the Friedan-Martinec-Shenker bosonization as a tool in studying commutant subalgebras of $\cS(V)$. I also thank A. Knutson, G. Schwarz, and N. Wallach for helpful discussions about classical invariant theory, especially the theory of invariant differential operators.

\newsec{Invariant differential operators}
Fix a basis $\{x_1,\dots,x_n\}$ for $V$ and a corresponding dual basis $\{x'_1,\dots,x'_n\}$ for $V^*$. The Weyl algebra $\cD(V)$ is generated by the linear functions $x_i'$ and the first-order differential operators $\frac{\partial}{\partial x'_i}$, which satisfy $[\frac{\partial}{\partial x'_i},x'_j] = \delta_{i,j}$. Equip $\cD(V)$ with the Bernstein filtration \eqn\bernstein{\cD(V)_{(0)}\subset \cD(V)_{(1)}\subset  \cdots,} defined by $(x'_1)^{k_1} \cdots (x'_n)^{k_n} (\frac{\partial}{\partial x'_1})^{l_1}\cdots (\frac{\partial}{\partial x'_n})^{l_n} \in \cD(V)_{(r)}$ if $k_1 + \cdots +k_n + l_1 + \cdots +l_n \leq r$. Given $\omega\in \cD(V)_{(r)}$ and $\nu\in\cD(V)_{(s)}$, $[\omega,\nu]\in \cD(V)_{(r+s-2)}$, so that \eqn\isodi{gr\cD(V) = \bigoplus_{r>0} \cD(V)_{(r)} / \cD(V)_{(r-1)} \cong Sym(V\oplus V^*).} We say that $deg(\alpha) = d$ if $\alpha\in\cD(V)_{(d)}$ and $\alpha\notin \cD(V)_{(d-1)}$. 

Let $G$ be a connected Lie group with Lie algebra $\gg$, and let $V$ be a linear representation of $G$ via $\rho: G\ra Aut(V)$. Then $G$ acts on $\cD(V)$ by algebra automorphisms, and induces an action $\rho^*: \gg \ra Der(\cD(V))$ by derivations of degree zero. Since $G$ is connected, the invariant ring $\cD(V)^G$ coincides with $\cD(V)^{\gg}$, where $$\cD(V)^{\gg} = \{\omega\in \cD(V)|~ \rho^*(\xi)(\omega) = 0,~ \forall \xi\in\gg\}.$$ We will usually work with the action of $\gg$ rather than $G$, and for greater flexibility, we do not assume that the $\gg$-action comes from an action of a {\it reductive} group $G$.

The action of $\gg$ on $\cD(V)$ can be realized by {\it inner} derivations: there is a Lie algebra homomorphism \eqn\defoftau{\tau:\gg\ra \cD(V),~~~~ \xi\mapsto - \sum_{i=1}^n x'_i \rho^*(\xi)\big(\frac{\partial}{\partial x'_i}\big) .} $\tau(\xi)$ is just the linear vector field on $V$ generated by $\xi$, so $\xi\in\gg$ acts on $\cD(V)$ by $[\tau(\xi),-]$. Clearly $\tau$ extends to a map $\gU\gg\ra \cD(V)$, and $$\cD(V)^{\gg}=Com(\tau(\gU\gg),\cD(V)).$$ Since $\gg$ acts on $\cD(V)$ by derivations of degree zero, \bernstein~restricts to a filtration
$\cD(V)^{\gg}_{(0)}\subset \cD(V)^{\gg}_{(1)} \subset \cdots$ on $\cD(V)^{\gg}$, and $gr(\cD(V)^{\gg}) \cong gr(\cD(V))^{\gg} \cong Sym(V\oplus V^*)^{\gg}$.

\subsec{The case where $\gg$ is abelian}

Our main focus is on the case where $\gg$ is the abelian Lie algebra $\C^m =  gl(1)\oplus\cdots\oplus gl(1)$, acting diagonally on $V$. Let $R(V)$ be the $\C$-vector space of all diagonal representations of $\gg$. Given $\rho\in R(V)$ and $\xi\in\gg$, $\rho(\xi)$ is a diagonal matrix with entries $a^{\xi}_1,\dots,a^{\xi}_n$, which we regard as a vector $a^{\xi} = (a^{\xi}_1,\dots, a^{\xi}_n)\in \C^n$. Let $A(\rho)\subset \C^n$ be the subspace spanned by $\{\rho(\xi)|~\xi\in\gg\}$. 

The action of $GL(m)$ on $\gg$ induces a natural action of $GL(m)$ on $R(V)$, defined by \eqn\actionofglm{(g\cdot \rho)(\xi) = \rho (g^{-1}\cdot \xi)} for all $g\in GL(m)$. Clearly $A(\rho) = A(g\cdot \rho)$ for all $g\in GL(m)$. Note that $dim~Ker(\rho) = dim~Ker (g\cdot \rho)$ for all $g\in GL(m)$, so in particular $GL(m)$ acts on the dense open set $R^0(V) = \{\rho\in R(V)|~ ker(\rho) = 0\}$. The correspondence $\rho\mapsto A(\rho)$ identifies $R^0(V)/GL(m)$ with the Grassmannian $Gr(m,n)$ of $m$-dimensional subspaces of $\C^n$.

Given $\rho\in R(V)$, $\cD(V)^{\gg} = \cD(V)^{\gg'}$ where $\gg' = \gg/Ker(\rho)$, so we may assume without loss of generality that $\rho\in R^0(V)$. We denote $\cD(V)^{\gg}$ by $\cD(V)^{\gg}_{\rho}$ when we need to emphasize the dependence on $\rho$. Given $\omega \in\cD(V)$, the condition $\rho^*(\xi)(\omega) = 0$ for all $\xi\in\gg$ is equivalent to the condition that $\rho^*( g\cdot \xi)(\omega) = 0$ for all $\xi\in\gg$, so it follows that $\cD(V)^{\gg}_{\rho} = \cD(V)^{\gg}_{g\cdot \rho}$ for all $g\in GL(m)$. Hence the family of algebras $\cD(V)^{\gg}_{\rho}$ is parametrized by the points $A(\rho)\in Gr(m,n)$.

Fix $\rho\in R^0(V)$, and choose a basis $\{\xi^1,\dots,\xi^m\}$ for $\gg$. Let  $a^i = (a^i_1,\dots,a^i_n)\in\C^n$ be the vectors corresponding to the diagonal matrices $\rho(\xi^i)$, and let $A = A(\rho)$ be the subspace spanned by these vectors. The map $\tau:\gg\ra \cD(V)$ is defined by
\eqn\taurho{\tau(\xi^i) = -\sum_{j=1}^n a^i_j x'_j \frac{\partial}{\partial x'_j}.} The Euler operators $\{e_j = x'_j \frac{\partial}{\partial x'_j}|~j=1,\dots,n\}$ lie in $\cD(V)^{\gg}$, and we denote the polynomial algebra $\C[e_1,\dots,e_n]$ by $E$.

For each $j=1,\dots,n$ and $d\in\Z$, define $v^d_j\in \cD(V)$ by
\eqn\hwvclassical{v_j^d = \bigg\{ \matrix{(\frac{\partial}{\partial x'_j})^{-d} & d<0 \cr 1 & d=0 \cr (x'_j)^d & d>0}.} Let $\Z^n\subset \C^n$ denote the lattice generated by the standard basis, and for each lattice point $l = (l_1,\dots,l_n) \in \Z^n$, define
\eqn\lmclassical{\omega_{l} = \prod_{j=1}^n v^{l_j}_j.} 
As a module over $E$, \eqn\modstructclassicali{\cD(V) = \bigoplus_{l\in\Z^n} M_l,} where $M_l$ is the free $E$-module generated by $\omega_l$. Moreover, we have
\eqn\modstructclassicalii{[e_j,\omega_l] = l_j \omega_l,} so the $\Z^n$-grading \modstructclassicali~is just the eigenspace decomposition of $\cD(V)$ under the family of diagonalizable operators $[e_j,-]$. In particular, \modstructclassicalii~shows that \eqn\perpcond{\rho^*(\xi^i)(\omega_l) = [\tau(\xi^i),\omega_l] = - \bra l,a^i\ket \omega_l,} where $\bra,\ket$ denotes the standard inner product on $\C^n$. Hence $\omega_l$ lies in $\cD(V)^{\gg}$ precisely when $l\in A^{\perp}$, so \eqn\modstructclassicaliii{\cD(V)^{\gg} = \bigoplus_{l\in A^{\perp}\cap \Z^n} M_l.} For generic actions, the lattice $A^{\perp}\cap \Z^n$ has rank zero, so $\cD(V)^{\gg} = M_0 = E.$

Consider the double commutant $Com(\cD(V)^{\gg},\cD(V))$, which always contains $T = \tau(\gU\gg) = \C[\tau(\xi_1)\dots,\tau(\xi_m)]$. Since $Com(E,\cD(V)) = E$, we have $Com(\cD(V)^{\gg},\cD(V)) = E$ for generic actions. 

Suppose next that $A^{\perp}\cap \Z^n$ has rank $r$ for some $0<r\leq n-m$. For $i=1,\dots,r$ let $\{l^i = (l^i_1,\dots,l^i_n)\}$ be a basis for $A^{\perp}\cap \Z^n$, and let $L$ be the $\C$-vector space spanned by $\{l^1,\dots,l^r\}$. If $r<n-m$, we can choose vectors $s^k = (s^k_1,\dots, s^k_n)\in L^{\perp} \cap A^{\perp}$, so that $\{l^1,\dots,l^r,s^{r+1},\dots,s^{n-m}\}$ is a basis for $A^{\perp}$. For $i=1,\dots,r$ and $k=r+1,\dots,n-m$, define differential operators $$\phi^i = \sum_{j=1}^n l^i_j e_j,~~~~~~~\psi^k =  \sum_{j=1}^n s^k_j e_j.$$ Note that $\C[e_1,\dots,e_n] = T\otimes \Psi\otimes \Phi$, where $\Phi = \C[\phi^1,\dots \phi^r]$ and $\Psi = \C[\psi^{r+1},\dots,\psi^{n-m}]$. 

\theorem{$Com(\cD(V)^{\gg},\cD(V)) = T\otimes \Psi$. Hence $\cD(V)^{\gg}$ and $T$ form a pair of mutual commutants inside $\cD(V)$ precisely when $\Psi = \C$, which occurs when $A^{\perp}\cap \Z^n$ has rank $n-m$.}\thmlab\howepairclassical

\proof By \modstructclassicalii, for any lattice point $l\in A^{\perp}\cap \Z^n$, and for $k=r+1,\dots,n-m$ we have $$[\psi^k,\omega_l] = \bra s^k,l\ket\omega_l = 0$$ since $s^k\in L^{\perp}$. It follows that $\Psi\subset Com(\cD(V)^{\gg},\cD(V))$. Hence $T\otimes \Psi\subset Com(\cD(V)^{\gg},\cD(V))$. Moreover, since $[\phi^i,\omega_l] = \bra l^i,l\ket\omega_l$ and $\{l^1,\dots,l^r\}$ form a basis for $A^{\perp}\cap \Z^n$, it follows that the variables $\phi^i$ cannot appear in any element $\omega\in Com(\cD(V)^{\gg},\cD(V))$. $\Box$

In the case $\Psi = \C$, we can recover the action $\rho$ (up to $GL(m)$-equivalence) from the algebra $\cD(V)^{\gg}$ by taking its commutant inside $\cD(V)$, but otherwise $\cD(V)^{\gg}$ does not determine the action.

\newsec{Vertex algebras}

We will assume that the reader is familiar with the basic notions in vertex algebra theory. For a list of references, see page 117 of \LL. We briefly describe the examples and constructions that we need, following the notation in \LL.

Given a Lie algebra $\gg$ equipped with a symmetric $\gg$-invariant bilinear form $B$, the {\it current algebra} $\cO(\gg,B)$ is the universal vertex algebra with generators $X^{\xi}(z)$, $\xi\in\gg$, which satisfy the OPE relations
$$X^{\xi}(z)X^{\eta}(w)\sim B(\xi,\eta)(z-w)^{-2} + X^{[\xi,\eta]}(w)(z-w)^{-1}.$$ 

Given a finite-dimensional vector space $V$, the $\beta\gamma$-system, or algebra of chiral differential operators $\cS(V)$, was introduced in \FMS. It is the unique vertex algebra with generators $\beta^{x}(z)$, $\gamma^{x'}(z)$ for $x\in V$, $x'\in V^*$, which satisfy
$$\beta^x(z)\gamma^{x'}(w)\sim\langle x',x\rangle (z-w)^{-1},~~~~~~~\gamma^{x'}(z)\beta^x(w)\sim -\langle x',x\rangle (z-w)^{-1},$$
\eqn\betagamma{\beta^x(z)\beta^y(w)\sim 0,~~~~~~~\gamma^{x'}(z)\gamma^{y'}(w)\sim 0.}
Given $\alpha = (\alpha_1,\dots,\alpha_n)\in\C^n$, $\cS(V)$ has a Virasoro element
\eqn\virbetagamma{L^{\alpha}(z) = \sum_{i=1}^n (\alpha_i-1):\partial\beta^{x_i}(z)\gamma^{x'_i}(z):\ +\alpha_i :\beta^{x_i}(z)\partial \gamma^{x'_i}(z):} of central charge $\sum_{i=1}^n \big(12\alpha_i^2 - 12\alpha_i + 2\big)$. Here $\{x_1,\dots, x_n\}$ is any basis for $V$ and $\{x'_1,\dots,x'_n\}$ is the corresponding dual basis for $V^*$. An OPE calculation shows that  $\beta^{x_i}(z),\gamma^{x'_i}(z)$ are primary of conformal weights $\alpha_i, 1-\alpha_i$, respectively. 

$\cS(V)$ has an additional $\Z$-grading which we call the $\beta\gamma$-charge. Define 
\eqn\defofv{v(z) = \sum_{i=1}^n : \beta^{x_i}(z)\gamma^{x'_i}(z):~.} The zeroth Fourier mode $v(0)$ acts diagonalizably on $\cS(V)$; the $\beta\gamma$-charge grading is just the eigenspace decomposition of $\cS(V)$ under $v(0)$. For $x \in V$ and $x'\in V^*$, $\beta^{x}(z)$ and $\gamma^{x'}(z)$ have $\beta\gamma$-charges $-1$ and $1$, respectively.

There is also an odd vertex algebra $\cE(V)$ known as a $bc$-system, or a semi-infinite
exterior algebra, which is generated by $b^x(z)$, $c^{x'}(z)$ for $x\in V$ and $x'\in V^*$, which satisfy
$$b^x(z)c^{x'}(w)\sim\langle x',x\rangle (z-w)^{-1},~~~~~~~c^{x'}(z) b^x(w)\sim \langle x',x\rangle (z-w)^{-1},$$ $$b^x(z)b^y(w)\sim 0,~~~~~~~c^{x'}(z)c^{y'}(w)\sim 0.$$
$\cE(V)$ has an analogous conformal structure $L^{\alpha}(z)$ for any $\alpha\in\C^n$, and an analogous $\Z$-grading which we call the $bc$-charge. Define \eqn\defofq{q(z) = -\sum_{i=1}^n : b^{x_i}(z) c^{x'_i}(z):~.} The zeroth Fourier mode $q(0)$ acts diagonalizably on $\cS(V)$, and the $bc$-charge grading is just the eigenspace decomposition of $\cE(V)$ under $q(0)$. Clearly $b^{x}(z)$ and $c^{x'}(z)$ have $bc$-charges $-1$ and $1$, respectively.

\subsec{The commutant construction}
\definition{Let $\cV$ be a vertex algebra, and let $\cA$ be a subalgebra. The commutant of $\cA$ in $\cV$, denoted by $Com(\cA,\cV)$ or $\cV^{\cA_+}$, is the subalgebra of vertex operators $v\in\cV$ such that $[a(z),v(w)] = 0$ for all $a\in\cA$. Equivalently, $a(z)\circ_n v(z) = 0$ for all $a\in\cA$ and $n\geq 0$.}
We regard $\cV$ as a module over $\cA$, and we regard $\cV^{\cA_+}$ as the invariant subalgebra. If $\cA$ is a homomorphic image of a current algebra $\cO(\gg,B)$, $\cV^{\cA_+}$ is just the invariant space $\cV^{\gg[t]}$. We will always assume that $\cV$ is equipped with a weight grading, and that $\cA$ is a graded subalgebra, so that $\cV^{\cA_+}$ is also a graded subalgebra of $\cV$.

Our main example of this construction comes from a representation $\rho: \gg\ra End(V)$ of a Lie algebra $\gg$. There is an induced vertex algebra homomorphism $\hat{\tau}:\cO(\gg,B)\ra \cS(V)$, which is analogous to the map $\tau: \gU \gg\ra \cD(V)$ given by \defoftau. Here $B$ is the bilinear form $B(\xi,\eta) = -Tr(\rho(\xi)\rho(\eta))$ on $\gg$. In terms of a basis $\{x_1,\dots, x_n\}$ for $V$ and dual basis $\{x'_1,\dots x'_n\}$ for $V^*$, $\hat{\tau}$ is defined by
\eqn\deftheta{\hat{\tau}(X^{\xi}(z)) = \theta^{\xi}(z) = - \sum_{i=1}^n~:\gamma^{x'_i}(z)\beta^{\rho(\xi)(x_i)}(z):~.} 

\definition{Let $\Theta$ denote the subalgebra $\hat{\tau}(\cO(\gg,B))\subset \cS(V)$. The commutant algebra $\cS(V)^{\Theta_+}$ will be called the algebra of {\it invariant chiral differential operators} on $V$.}

If $\cS(V)$ is equipped with the conformal structure $L^{\alpha}$ given by \virbetagamma, $\Theta$ is not a graded subalgebra of $\cS(V)$ in general. For example, if $\gg = gl(n)$ and $V=\C^n$, $\Theta$ is graded by weight precisely when $\alpha_1 = \alpha_2 = \cdots  = \alpha_n$. However, when $\gg$ is abelian and its action on $V$ is diagonal, $\theta^{\xi}(z)$ will be homogeneous of weight one for any $\alpha$. Hence $\cS(V)^{\Theta_+}$ is also graded by weight, but this grading will depend on the choice of $\alpha$.

\subsec{The Zhu functor}
Let $\cV$ be a vertex algebra with weight grading $\cV = \bigoplus_{n\in\Z} \cV_n$. In \Zh, Zhu introduced a functor that attaches to $\cV$ an associative algebra $A(\cV)$, together with a surjective linear map $\pi_{Zh}:\cV\ra A(\cV)$. For $a\in \cV_{m}$ and $b\in\cV$, we define
\eqn\defzhu{a*b = Res_z \bigg (a(z) \frac{(z+1)^{m}}{z}b\bigg),} and extend $*$ by linearity to a bilinear operation $\cV\otimes \cV\ra \cV$. Let $O(\cV)$ denote the subspace of $\cV$ spanned by elements of the form
\eqn\zhuideal{a\circ b = Res_z \bigg (a(z) \frac{(z+1)^{m}}{z^2}b\bigg)} where $a\in \cV_m$, and let $A(\cV)$ be the quotient $\cV/O(\cV)$, with projection $\pi_{Zh}:\cV\ra A(\cV)$. For $a,b\in \cV$, $a\sim b$ means $a-b\in O(\cV)$, and $[a]$ denotes the image of $a$ in $A(\cV)$. A useful fact which is immediate from \defzhu~and \zhuideal~is that for $a\in\cV_m$, \eqn\uf{\partial a \sim ma.}

\theorem{(Zhu) $O(\cV)$ is a two-sided ideal in $V$ under the product $*$, and $(A(\cV),*)$ is an associative algebra with unit $[1]$. The assignment $\cV\mapsto A(\cV)$ is functorial. If $\cI$ is a vertex algebra ideal of $\cV$, we have \eqn\zhuidealrule{A(\cV/\cI)\cong A(\cV)/ I,~~~~~I = \pi_{Zh}(\cI).}}

The main application of the Zhu functor is to study the representation theory of $\cV$, or at least reduce it to a more classical problem. Let $M = \bigoplus_{n\geq 0} M_n$ be a module over $\cV$ such that for $a\in\cV_m$, $a(n) M_k \subset M_{m+k -n-1}$ for all $n\in\Z$. Given $a\in\cV_m$, the Fourier mode $a(m-1)$ acts on each $M_k$. The subspace $M_0$ is then a module over $A(\cV)$ with action $[a]\mapsto a(m-1) \in End(M_0)$. In fact, $M\mapsto M_0$ provides a one-to-one correspondence between irreducible $\Z_{\geq 0}$-graded $\cV$-modules and irreducible $A(\cV)$-modules. 

A vertex algebra $\cV$ is said to be {\it strongly generated} by a subset $\{v_i(z)|~i\in I\}$ if $\cV$ is spanned by collection of iterated Wick products $$\{ :\partial^{k_1} v_{i_1}(z)\cdots \partial^{k_m} v_{i_m}(z):~|~ k_1,\dots,k_m \geq 0\}.$$
\lemma{Suppose that $\cV$ is strongly generated by $\{v_i(z)|~i\in I\}$, which are homogeneous of weights $d_i\geq 0$. Then $A(\cV)$ is generated as an associative algebra by the collection $\{\pi_{Zh}(v_i)|~i\in I\}$.}\thmlab\zhugeneration

\proof Let $\cC$ be the algebra generated by $\{\pi_{Zh}(v_i)| i\in I\}$. We need to show that for any vertex operator $\omega\in \cV$, we have $\pi_{Zh}(\omega)\in \cC$. By strong generation, it suffices to prove this when $\omega$ is a monomial of the form $$:\partial^{k_1} v_{i_1}\cdots \partial^{k_r} v_{i_r}:~.$$ We proceed by induction on weight. Suppose first that $\omega$ has weight zero, so that $k_1 = \cdots = k_r = 0$ and $v_{i_1},\dots,v_{i_r}$ all have weight zero. Note that $v_{i_1}\circ_n (:v_{i_2}\cdots v_{i_r}:)$ has weight $-n-1$, and hence vanishes for all $n\geq 0$. It follows from \defzhu~that $$[v_{i_1}]* [:v_{i_2}\cdots v_{i_r}:] = [\omega].$$ Continuing in this way, we see that $[\omega] = [v_{i_1}]*[ v_{i_2}]*\cdots * [v_{i_r}]\in\cC$. Next, assume that $\pi_{Zh}(\omega)\in\cC$ whenever $wt(\omega)<n$, and suppose that $\omega =  ~:\partial^{k_1} v_i \cdots \partial^{k_r} v_r:~$ has weight $n$. We calculate
$$[\partial^{k_1} v_{i_1}]* [ : \partial^{k_2}v_{i_2}\cdots \partial^{k_r} v_{i_r}:] = [\omega] + \cdots,$$ where $\cdots$ is a linear combination of terms of the form $[\partial^{k_1} v_{i_1} \circ_k (: \partial^{k_2}v_{i_2}\cdots \partial^{k_r} v_{i_r}:)]$ for $k\geq 0$. The vertex operators $\partial^{k_1} v_{i_1} \circ_k (: \partial^{k_2}v_{i_2}\cdots \partial^{k_r} v_{i_r}:)$ all have weight $n-k-1$, so by our inductive assumption, $[\partial^{k_1} v_{i_1} \circ_k (: \partial^{k_2}v_{i_2}\cdots \partial^{k_r} v_{i_r}:)]\in\cC$. Applying the same argument to the vertex operator $: \partial^{k_2}v_{i_2}\cdots \partial^{k_r} v_{i_r}:$ and proceeding by induction on $r$, we see that $[\omega] \equiv [\partial^{k_1} v_{i_1}]*\cdots *[\partial^{k_n} v_{i_n}]$ modulo $\cC$. Finally, by applying \uf~repeatedly, we see that $[\omega]\in\cC$, as claimed. $\Box$.

\others{Example}{$\cV = \cO(\gg,B)$ where each generator $X^{\xi}$ has weight 1. Then $A(\cO(\gg,B))$ is generated by $\{[X^{\xi}]|~\xi\in\gg\}$, and is isomorphic to the universal enveloping algebra $\gU\gg$ via $[X^{\xi}]\mapsto \xi$.}

\others{Example}{Let $\cV = \cS(V)$ where $V=\C^n$, and $\cS(V)$ is equipped with the conformal structure $L^{\alpha}$ given by \virbetagamma. Then $A(\cS(V))$ is generated by $\{ [\gamma^{x'_i}], [\beta^{x_i}] \}$ and is isomorphic to the Weyl algebra $\cD(V)$ with generators $x'_i, \frac{\partial}{\partial x'_i}$ via $$[\gamma^{x'_i}]\mapsto x'_i,~~~~ [\beta^{x_i}]\mapsto \frac{\partial}{\partial x'_i}.$$
Even though the structure of $A(\cS(V))$ is independent of the choice of $\alpha$, the Zhu map $\pi_{Zh}:\cS(V)\ra A(\cS(V))$ does depend on $\alpha$. For example, \defzhu~shows that \eqn\pizhubg{\pi_{Zh}(:\gamma^{x'_i}\beta^{x_i}:) =  x'_i\frac{\partial}{\partial x'_i} +1- \alpha_i.}}

We will be particularly concerned with the interaction between the commutant construction and the Zhu functor. If $a,b\in\cV$ are (super)commuting vertex operators, $[a]$ and $[b]$ are (super)commuting elements of $A(\cV)$. Hence for any subalgebra $\cB\subset \cV$, we have a commutative diagram
\eqn\cdgencase{\matrix{Com(\cB,\cV)&\br\iota\over\hookrightarrow&\cV \cr \pi\da\hskip.4in& &\pi_{Zh}\da\hskip.4in\cr Com(B,A(\cV))&\br\iota\over\hookrightarrow & A(\cV)}.} Here $B$ denotes the subalgebra $\pi_{Zh}(\cB)\subset A(\cV)$, and $Com(B,A(\cV))$ denotes the (super)commutant of $B$ inside $A(\cV)$. The horizontal maps are inclusions, and $\pi$ is the restriction of the Zhu map on $\cV$ to $Com(\cB,\cV)$. Clearly $Im(\pi)$ is a subalgebra of $Com(B,A(\cV))$. A natural problem is to describe $Im(\pi)$ and $Coker(\pi)$. In our main example $\cV = \cS(V)$ and $\cA= \Theta$, we have $\pi_{Zh}(\Theta)= \tau(\gU\gg)\subset \cD(V)$ and $Com(\tau(\gU\gg),\cD(V)) = \cD(V)^{\gg}$, so \cdgencase~specializes to \commdiagva.

\newsec{The Friedan-Martinec-Shenker bosonization}

\subsec{Bosonization of fermions}
First we describe the bosonization of fermions and the well-known boson-fermion correspondence due to \FI. Let $A$ be the Heisenberg algebra with generators $j(n)$, $n\in \Z$, and $\kappa$, satisfying $[j(n),j(m)] = n \delta_{n+m,0}\kappa$. The field $j(z)= \sum_{n\in \Z} j(n) z^{-n-1}$ satisfies the OPE $$j(z)j(w)\sim (z-w)^{-2},$$ and generates a Heisenberg vertex algebra $\cH$ of central charge $1$. Define the {\it free bosonic scalar field}
$$\phi(z) = q+ j(0) \ln z - \sum_{n\neq 0} \frac{j(n)}{n} x^{-n},$$ where $q$ satisfies $[q,j(n)] = \delta_{n,0}$. Clearly $\partial \phi(z) = j(z)$, and we have the OPE $$\phi(z)\phi(w)\sim \ln(z-w).$$
Given $\alpha\in\C$, let $\cH_{\alpha}$ denote the irreducible representation of $A$ generated by the vacuum vector $v_{\alpha}$ satisfying \eqn\heismodule{j(n) v_{\alpha}= \alpha \delta_{n,0} v_{\alpha},~~~ n\geq 0.} Given $\eta\in\C$, the operator $e^{\eta q}(v_{\alpha}) = v_{\alpha+\eta}$, so $e^{\eta q}$ maps $\cH_{\alpha}\ra \cH_{\alpha+\eta}$. Define the vertex operator
$$X_{\eta}(z) = e^{\eta \phi(z)} = e^{\eta q} z^{\eta \alpha} exp(\eta \sum_{n>0} j(-n) \frac{z^n}{n}) exp(\eta\sum_{n<0} j(-n) \frac{z^n}{n}).$$ The $X_{\eta}$ satisfy the OPEs
$$j(z) X_{\eta}(w) = \eta X_{\eta}(w)(z-w)^{-1} + \frac{1}{\eta} \partial X_{\eta}(w),$$
$$ X_{\eta}(z) X_{\nu}(w) = (z-w)^{\eta\nu} :X_{\eta}(z) X_{\nu}(w):.$$
If we take $\eta = \pm 1$, the pair of (fermionic) fields $X_1,X_{-1}$ generate the lattice vertex algebra $V_L$ associated to the one-dimensional lattice $L=\Z$. The state space of $V_L$ is just $\sum_{n\in\Z} \cH_n = \cH\otimes_{\C} L$. It follows that $$X_1(z)X_{-1}(w) \sim (z-w)^{-1},~~~~~~~X_{-1}(z)X_{1}(w) \sim (z-w)^{-1},$$ $$X_1(z)X_1(w) \sim 0,~~~~~~~ X_{-1}(z)X_{-1}(w)\sim 0,$$  so the map $\cE\ra V_L$ sending $b\mapsto X_{-1}, c\ra X_1$ is a vertex algebra isomorphism. Here $\cE$ denotes the $bc$-system $\cE(V)$ in the case where $V$ is one-dimensional.

\subsec{Bosonization of bosons}
Next, we describe the bosonization of bosons, following \FF. Recall that $\cE$ has the grading $\cE = \oplus_{l\in\Z} \cE^l$ by $bc$-charge. As in \FF, define $N(s) = \sum_{l\in\Z} \cE^l\otimes \cH_{i(s+l)}$, which is a module over the vertex algebra $\cE\otimes V_{L'}$. Here $L'$ is the one-dimensional lattice $i\Z$, and $V_{L'}$ is generated by $X_{\pm i}$. We define a map $\epsilon:\cS\ra \cE\otimes V_{L'}$ by 
\eqn\eps{\beta \mapsto \partial b \otimes X_{-i},~~~~~~~\gamma \mapsto  c\otimes X_i.}
It is straightforward to check that \eps~is a vertex algebra homomorphism, which is injective since $\cS$ is simple. Moreover Proposition 3 of \FF~shows that the image of \eps~coincides with the kernel of $c(0): N(s)\ra N(s-1)$. Let $\cE'$ be the subalgebra of $\cE$ generated by $c$ and $\partial b$, which coincides with the kernel of $c(0):\cE\ra \cE$. It follows that \eqn\imageofeps{\epsilon(\cS)\subset \cE'\otimes V_{L'}.}

\newsec{$\cW$ algebras}

The $\cW$ algebras are vertex algebras which arise as extended symmetry algebras of two-dimensional conformal field theories. For each integer $n\geq 2$ and $c\in\C$, the algebra $\cW_{n,c}$ of central charge $c$ is generated by fields of conformal weights $2,3,\dots,n$. In the case $n=2$, $\cW_{2,c}$ is just the Virasoro algebra of central charge $c$. In contrast to the Virasoro algebra, the generating fields for $\cW_{n,c}$ for $n\geq 3$ have nonlinear terms in their OPEs, which makes the representation theory of these algebras highly nontrivial. One also considers various limits of $\cW$ algebras denoted by $\cW_{1+\infty,c}$ which may be defined as modules over the universal central extension $\hat{\cD}$ of the Lie algebra $\cD$ of differential operators on the circle \KR. 

We will be particularly concerned with the $\cW_3$ algebra, which was introduced by Zamolodchikov in \Za~and studied extensively in \BMP. Our discussion is taken directly from \WI\WII. First, let $\cF(\cW_3)$ denote the free associative algebra with generators $L_m,W_m$, $m\in\Z$. Let $\hat{\cF}(\cW_3)$ be the completion of $\cF(\cW_3)$ consisting of (possibly) infinite sums of monomials in $\cF(\cW_3)$ such that for each $N>0$, only finitely many terms depend only on the variables $L_n,W_n$ for $n\leq N$. For a fixed central charge $c\in\C$, let $\gU\cW_{3,c}$ be the quotient of $\hat{\cF}(\cW_3)$ by the ideal generated by
\eqn\wcommi{[L_m,L_n] = (m-n)L_{m+n} + \frac{c}{12} (m^3-m) \delta_{m,-n},}
\eqn\wcommii{[L_m,W_n] = (2m-n) W_{m+n},}
\eqn\wcommiii{[W_m,W_n] = (m-n) \bigg( \frac{1}{15} (m+n+3)(m+n+2) -\frac{1}{6}(m+2)(n+2)\bigg) L_{m+n}} $$+ \frac{16}{22+5c}(m-n) \Lambda_{m+n} + \frac{c}{360} m(m^2-1)(m^2-4) \delta_{m,-n}.$$ Here $$\Lambda_m = \sum_{n\leq -2} L_n L_{m-n} + \sum_{n>-2} L_{m-n} L_n - \frac{3}{10} (m+2)(m+3) L_m.$$ Let $$\cW_{3,c,\pm} = \{L_n,W_n,~ \pm n >0\},~~~~ \cW_{3,c,0} = \{L_0,W_0\}.$$ The Verma module $\cM_c(t,w)$ of highest weight $(t,w)$ is the induced module $$\gU\cW_{3,c}\otimes_{\cW_{3,c,+}\oplus \cW_{3,c,0}} \C_{t,w},$$ where $\C_{t,w}$ is the one-dimensional $\cW_{3,c,+}\oplus \cW_{3,c,0}$-module generated by the vector $v_{t,w}$ such that
$$\cW_{3,c,+}(v_{t,w}) = 0,~~~ L_0(v_{t,w}) = t v_{t,w},~~~ W_0 (v_{t,w}) = w v_{t,w}.$$ A vector $v\in \cM_c(t,w)$ is called {\it singular} if $\cW_{3,c,+}(v) = 0$. In the case $t = w = 0$, the vectors \eqn\singvect{L_{-1}(v_{0,0}),~~~ W_{-1}(v_{0,0}),~~~W_{-2}(v_{0,0})} are singular vectors in $\cM_c(0,0)$. The {\it vacuum module} $\cV\cW_{3,c}$ is defined to be the quotient of $\cM_c(0,0)$ by the $\gU\cW_{3,c}$-submodule generated by the vectors \singvect. $\cV\cW_{3,c}$ has the structure of a vertex algebra which is freely generated by the vertex operators
$$L(z) = \sum_{n\in\Z} L_n z^{-n-2},~~~ W(z) = \sum_{n\in\Z} W_n z^{-n-3}.$$ In particular, the vertex operators $$\{\partial^{i_1} L(z)\cdots \partial^{i_m} L(z) \partial^{j_1} W(z)\cdots \partial^{j_n} W(z)|~ 0\leq i_1\leq\cdots\leq i_m,~~ 0\leq j_1\leq\cdots\leq j_n\}$$ which correspond to 
$ i_1!\cdots i_m! j_1!\cdots j_n!L_{-i_1 - 2}\cdots L_{-i_m - 2} W_{-j_1 -3}\cdots W_{-j_n -3}v_{0,0}$ under the state-operator correspondence, form a basis for $\cV\cW_{3,c}$. By Lemma 4.1 of \WII, the Zhu algebra $A(\cV\cW_{3,c})$ is just the polynomial algebra $\C[l,w]$ where $l = \pi_{Zh}(L)$ and $w = \pi_{Zh}(W)$. 

Let $\cI_c$ denote the maximal proper $\gU\cW_{3,c}$-submodule of $\cV\cW_{3,c}$, which is a vertex algebra ideal. The quotient $\cV\cW_{3,c}/\cI_c$ is a simple vertex algebra which we denote by $\cW_{3,c}$. Let $I_c = \pi_{Zh}(\cI_c)$, which is an ideal of $\C[l,w]$. By \zhuidealrule, we have $A(\cW_{3,c}) = \C[l,w] / I_c$. Generically, $\cI_c = 0$, so that $\cV\cW_{3,c} = \cW_{3,c}$. We will be primarily concerned with the non-generic case $c=-2$, in which $\cI_{-2}\neq 0$. 
The generators $L(z),W(z)\in \cV\cW_{3,-2}$ satisfy the following OPEs:
\eqn\wi{L(z) L(w) \sim -(z-w)^{-4} + 2 L(w)(z-w)^{-2} + \partial L(w) (z-w)^{-1},}
\eqn\wii{L(z) W(w)\sim 3 W(w)(z-w)^{-2} + \partial W(w)(z-w)^{-1},}
$$W(z)W(w)\sim -\frac{2}{3} (z-w)^{-6} + 2L(w) (z-w)^{-4} + \partial L(w)(z-w)^{-3}$$ $$+ (\frac{8}{3} :L(w)L(w): - \frac{1}{2} \partial^2 L(w))(z-w)^{-2} $$ \eqn\wiii{+(\frac{4}{3} \partial(:L(w)L(w):) - \frac{1}{3} \partial^3 L(w))(z-w)^{-1}.} The simple vertex algebra $\cW_{3,-2}$ also has generators $L(z),W(z)$ satisfying \wi-\wiii, but $\cW_{3,-2}$ is no longer freely generated. 

In order to avoid introducing extra notation, we will {\it not} use the change of variables $\tilde{W}(z) = \frac{1}{2}\sqrt{6} W(z)$ given by Equation 3.13 of \WII. By Lemma 4.3 of \WII, the ideal $I_{-2}\subset \C[l,w]$ is generated (in our variables) by the polynomial 
\eqn\zwideal{w^2 - \frac{2}{27} l^2(8l+1).}

\subsec{The representation theory of $\cW_{3,-2}$}
In \WII, W. Wang gave a complete classification of the irreducible modules over the simple vertex algebra $\cW_{3,-2}$. An important ingredient in his classification is the following realization of $\cW_{3,-2}$ as a subalgebra of the Heisenberg algebra $\cH$ with generator $j(z)$ satisfying $j(z)j(w)\sim (z-w)^{-2}$. Define 
\eqn\heisi{L_{\cH} = \frac{1}{2}(: j^2:) + \partial j,~~~~~~~~~W_{\cH} =\frac{2}{3\sqrt{6}}(: j^3:) ~+~ \frac{1}{\sqrt{6}} (:j \partial j:)~+~ \frac{1}{6\sqrt{6}} \partial^2 j.}
The map $\cW_{3,-2}\hookrightarrow \cH$ sending $L\mapsto L_{\cH}$ and $W\mapsto W_{\cH}$ is a vertex algebra homomorphism, so we may regard any $\cH$-module as a $\cW_{3,-2}$-module. Given $\alpha\in\C$, consider the irreducible $\cH$-module $\cH_{\alpha}$ defined by \heismodule, and let $V_{\alpha}$ denote the irreducible quotient of the $\cW_{3,-2}$-submodule of $\cH_{\alpha}$ generated by $v_{\alpha}$. It is easily checked that the generator $v_{\alpha}$ is a highest weight vector of $\cW_{3,-2}$ with highest weight 
\eqn\highweight{\bigg(\frac{1}{2} \alpha (\alpha -1), \frac{1}{3\sqrt{6}} \alpha (\alpha -1)(2\alpha -1)\bigg).}
The main result of \WII~is that the modules $\{V_{\alpha}|~\alpha\in\C\}$ account for all the irreducible modules of $\cW_{3,-2}$.

\newsec{The commutant algebra $\cS(V)^{\Theta_+}$ for $\gg = gl(1)$ and $V = \C$}

In this section, we describe $\cS(V)^{\Theta_+}$ in the case where $\gg = gl(1)$ and $V = \C$, where the action $\rho:\gg\ra End ~V$ is by multiplication. Fix a basis $\xi$ of $\gg$ and a basis $x$ of $V$, such that $\rho(\xi)(x) = x$. Then $\cS= \cS(V)$ is generated by $\beta(z) = \beta^x(z)$ and $\gamma(z) = \gamma^{x'}(z)$, and the map \taurho~is given by $$\gg\ra \cD = \cD(V),~~~~~~~\xi\mapsto -x' \frac{d}{dx'}.$$ In this case, $\cO(\gg,B)$ is just the Heisenberg algebra $\cH$ of central charge $-1$, and the action of $\cH$ on $\cS$ given by \deftheta~is 
\eqn\thetaonedimcase{ \theta(z) = -:\gamma(z)\beta(z):~,} which clearly satisfies \eqn\thetatheta{\theta(z)\theta(w)\sim -(z-w)^{-2}.}  As usual, $\Theta$ will denote the subalgebra of $\cS$ generated by $\theta(z)$. Since $-\theta(0)$ is the $\beta\gamma$-charge operator, $\cS^{\Theta_+}$ must lie in the subalgebra $\cS^0$ of $\beta\gamma$-charge zero.

Let $:\theta^n:$ denote the $n$-fold iterated Wick product of $\theta$ with itself. It is clear from \thetatheta~that each $:\theta^n:$ lies in $\cS^0$ but not in $\cS^{\Theta_+}$. A natural place to look for elements in $\cS^{\Theta_+}$ is to begin with the operators $:\theta^n:$ and try to \lq\lq quantum correct" them so that they lie in $\cS^{\Theta_+}$. As a polynomial in $\beta,\partial\beta,\dots,\gamma,\partial\gamma,\cdots$, note that $$:\theta^n:~ = (-1)^n\beta^n\gamma^n + \nu_n,$$ where $\nu_n$ has degree at most $2n-2$. By a quantum correction, we mean an element $\omega_n\in \cS$ of polynomial degree at most $2n-2$, so that $:\theta^n: + ~\omega_n\in \cS^{\Theta_+}$.  

Clearly $\theta$ has no such correction $\omega_1$, because $\omega_1$ would have to be a scalar, in which case $\theta\circ_1(\theta + \omega_1) = \theta\circ_1\theta = -1$. However, the next lemma shows that we can find such $\omega_n$ for all $n\geq 2$.

\lemma{Let $$\omega_2 = ~: \beta(\partial\gamma):~  - ~:(\partial\beta)\gamma:~,$$ $$\omega_3 = -\frac{9}{2}~:\beta^2 \gamma(\partial\gamma): ~+~\frac{9}{2} :\beta(\partial\beta)\gamma^2: ~-~ \frac{3}{2}:\beta(\partial^2\gamma): ~-~ \frac{3}{2}:(\partial^2\beta)\gamma:~ +6  ~:(\partial \beta)(\partial\gamma):~.$$
Then $:\theta^2:~ +~ \omega_2\in \cS^{\Theta_+}$ and $:\theta^3:~+~\omega_3\in\cS^{\Theta_+}$. Since $:(\theta^n):$ and $:(:\theta^i:)(:\theta^j:):$ have the same leading term as polynomials in $\beta,\partial\beta,\dots,\gamma,\partial\gamma,\cdots$ for $i+j=n$, it follows that for any $n\geq 2$ we can find $\omega_n$ such that $:\theta^n:+~\omega_n\in\cS^{\Theta_+}$.}\thmlab\qcor

\proof This is a straightforward OPE calculation. $\Box$
 
Next, define vertex operators $L_{\cS},W_{\cS}\in\cS^{\Theta_+}$ as follows:
\eqn\LS{L_{\cS} = \frac{1}{2}(:\theta^2:~ + ~\omega_2) =  \frac{1}{2}(:\beta^2\gamma^2:)~-~:(\partial \beta) \gamma: ~+~ :\beta (\partial \gamma):~,}
$$W_{\cS} = - \sqrt{\frac{2}{27}} (:\theta^3:~ +~ \omega_3) $$ $$ = \sqrt{\frac{2}{27}}(:\beta^3\gamma^3:)~ -~ \sqrt{\frac{3}{2}}(: \beta(\partial\beta)\gamma^2:)~ + ~\sqrt{\frac{3}{2}}( : \beta^2\gamma(\partial \gamma):) $$ \eqn\WS{+ ~\sqrt{\frac{1}{6}}(:(\partial^2\beta)\gamma:)~ -~\sqrt{\frac{8}{3}}(: (\partial\beta)(\partial\gamma):)~ +~\sqrt{\frac{1}{6}}( :\beta(\partial^2\gamma):).} 
Let $\cW\subset \cS^{\Theta_+}$ be the vertex algebra generated by $L_{\cS},W_{\cS}$. An OPE calculation shows that the map \eqn\deff{\cV\cW_{3,-2}\ra \cS^{\Theta_+},~~~~ L\mapsto L_{\cS},~~~~W\mapsto W_{\cS}} is a vertex algebra homomorphism. Moreover, the ideal $\cI_{-2}$ is annihilated by \deff, so this map descends to a map \eqn\defff{f:\cW_{3,-2}\hookrightarrow \cS^{\Theta_+}.}
In fact, \defff~is related to the realization of $\cW_{3,-2}$ as a subalgebra of $\cH$ defined earlier. First, under the boson-fermion correspondence, \eqn\bfci{L_{\cH}\mapsto L_{\cE} = ~:\partial b c:~,}
\eqn\bfcii{W_{\cH} \mapsto W_{\cE} = \frac{1}{\sqrt{6}}( : (\partial^2 b) c:~-~:(\partial b)( \partial c):).} Next,
under the map $\epsilon: \cS\ra \cE\otimes \cH$ given by \eps, we have \eqn\imepsilon{L_{\cS}\mapsto L_{\cE}\otimes 1,~~~~W_{\cS}\mapsto W_{\cE}\otimes 1.}

The subalgebra $\cS^0$ of $\beta\gamma$-charge zero has a natural set of generators $$\{J^i = ~:\beta(\partial^i\gamma):~, i\geq 0\},$$ and it is well known that $\cS^0$ is isomorphic to $\cW_{1+\infty, -1}$ \KR. One of the main results of \WI~is that $\epsilon:\cS\ra \cE\otimes \cH$ restricts to an isomorphism \eqn\wang{\cS^0\cong \cA\otimes \cH,} where $\cA\cong \cW_{3,-2}$ is the subalgebra of $\cE$ generated by $L_{\cE}$ and $W_{\cE}$. By \imepsilon, $\epsilon$ maps $\cW$ onto $\cA\otimes 1$. Similarly, $\epsilon (\theta) = i (1\otimes j)$, so $\epsilon$ maps $\Theta$ onto $1\otimes \cH$, and $\cS^0 = \cW\otimes \Theta$.

For each $d\in\Z$, the subspace $\cS^d$ of $\beta\gamma$-charge $d$ is a module over $\cS^0$, which is in fact irreducible \KR\WII. Define $v^d(z)\in \cS^d$ by
\eqn\highweightvect{v^d(z) = \bigg\{ \matrix{\beta(z)^{-d} & d<0 \cr 1 & d=0 \cr \gamma(z)^d & d>0}.}
Here $\beta(z)^{-d}$ and $\gamma(z)^d$ denote the $d$-fold iterated Wick products $:\beta(z)\cdots\beta(z):$ and $:\gamma(z)\cdots\gamma(z):$, respectively. Each $v^d(z)$ is a highest weight vector for the action of $\cW_{3,-2}$, and the highest weight of $v^d(z)$ is given by \highweight~with
\eqn\weightsofvd{\bigg\{ \matrix{\alpha = d & d\leq0 \cr \alpha = d+1 & d>0}.} Moreover, $v^d(z)$ is also a highest weight vector for the action of $\cH$, so $\cS^d$ is generated by $v^d(z)$ as a module over $\cW_{3,-2}\otimes \cH$.

\theorem{The map $f:\cW_{3,-2}\hookrightarrow \cS^{\Theta_+}$ given by \defff~is an isomorphism of vertex algebras. Moreover, $Com(\cS^{\Theta_+},\cS) = \Theta$. Hence $\Theta$ and $\cS^{\Theta_+}$ form a Howe pair inside $\cS$.}\thmlab\main

\proof Clearly $\cS^{\Theta_+}\subset \cS^0$, and since $\cS^0 = \cW\otimes \Theta$, we have $$\cS^{\Theta_+} = Com(\Theta, \cW\otimes \Theta) = \cW\otimes Com(\Theta,\Theta) = \cW.$$ This proves the first statement. As for the second statement, it is clear from \highweight~and \weightsofvd~that $Com(\cS^{\Theta_+},\cS)\subset \cS^0$. Hence $$Com(\cS^{\Theta_+},\cS) = Com (\cW, \cW\otimes \Theta) = \Theta\otimes Com(\cW,\cW) = \Theta.~~~\Box$$

\subsec{The map $\pi: \cS^{\Theta_+}\ra \cD^{\gg}$}
Equip $\cS$ with the conformal structure $L^{\alpha} = (\alpha-1) :\partial \beta(z) \gamma(z):~ + \alpha :\beta(z) \partial \gamma(z):$, and consider the map $\pi:\cS^{\Theta_+}\ra \cD^{\gg}$ given by \commdiagva. In this case, $\cD^{\gg}$ is just the polynomial algebra $\C[e]$, where $e$ is the Euler operator $x' \frac{d}{dx'}$. 
\lemma{We have \eqn\pzhu{\pi(L_{\cS}) = \frac{1}{2}(e^2 + e),~~~~~~~ \pi(W_{\cS}) = \frac{2}{3\sqrt{6}} e^3 +\frac{1}{\sqrt{6}} e^2 + \frac{1}{3\sqrt{6}} e.} In particular, $\pi(L_{\cS})$ and $\pi(W_{\cS})$ are independent of the choice of $\alpha$.}\thmlab\zhugeni 

\proof This is a straightforward computation using \defzhu~and the fact that $\pi_{Zh}(\gamma(z)) = x'$ and $\pi_{Zh}(\beta(z)) = \frac{d}{dx'}$. Note that $l=\pi(L_{\cS})$ and $w=\pi(W_{\cS})$ satisfy \zwideal. $\Box$

\corollary{For any conformal structure $L^{\alpha}$ on $\cS$ as above, $Im(\pi)$ is the subalgebra of $\C[e]$ generated by $\pi(L_{\cS})$ and $\pi(W_{\cS})$. Moreover, $Coker(\pi) = \C[e]/ Im(\pi)$ has dimension one, and is spanned by the image of $e$ in $Coker(\pi)$.}

\proof The first statement is immediate from Lemma \zhugeneration, since $\cS^{\Theta_+}$ is strongly generated by $L_{\cS}$ and $W_{\cS}$ which have weights $2$ and $3$ respectively. The second statement follows from \pizhubg~and \pzhu, because any polynomial in $\C[e]$ is equivalent to an element which is homogeneous of degree 1 modulo $Im(\pi)$. $\Box$

\newsec{$\cS(V)^{\Theta_+}$ for abelian Lie algebra actions}

Fix a basis $\{x_1,\dots,x_n\}$ for $V$ and dual basis $\{x'_1,\dots,x'_n\}$ for $V^*$. We regard $\cS(V)$ as $\cS_1\otimes \cdots\otimes \cS_n$, where $\cS_j$ is the copy of $\cS$ generated by $\beta^{x_j}(z),\gamma^{x'_j}(z)$. Let $f_j:\cS\ra \cS(V)$ be the obvious map onto the $j$th factor. 
The subspace $\cS_j^0$ of $\beta\gamma$-charge zero is isomorphic to $\cW^j\otimes \cH^j$, where $\cH^j$ is generated by $\theta^j(z) = f_j(\theta(z))$, and $\cW^j$ is generated by $L^j = f_j (L_{\cS})$, $W^j = f_j(W_{\cS})$. Moreover, as a module over $\cW^j\otimes \cH^j$, the space $\cS_j^d$ of $\beta\gamma$-charge $d$ is generated by the highest weight vector $v^d_j(z)= f_j(v^d(z))$, which is given by \eqn\hwv{v_j^d(z) = \bigg\{ \matrix{\beta^{x_j}(z)^{-d} & d<0 \cr 1 & d=0 \cr \gamma^{x'_j}(z)^d & d>0}.}We denote by $\cS'_j$ the linear span of the vectors $\{v^d_j(z)|~d\in\Z\}$. Note that for any conformal structure $L^{\alpha}$ on $\cS(V)$, the differential operators $v^d_j\in\cD(V)$ defined by \hwvclassical~correspond to $v^d_j(z)$ under the Zhu map. Let $\cB$ denote the vertex algebra $$\cS_1^0\otimes \cdots\otimes\cS_n^0 \cong (\cW^1\otimes \cH^1)\otimes \cdots \otimes (\cW^n\otimes \cH^n).$$ Clearly the space $\cS(V)'$ consisting of highest-weight vectors for the action of $\cB$ is just $\cS'_1\otimes \cdots\otimes \cS'_n$. As usual, let $\Z^n\subset \C^n$ denote the standard lattice. For each lattice point $l = (l_1,\dots,l_n) \in \Z^n$, define
\eqn\latticemap{\omega_{l}(z) = ~: v^{l_1}_1(z)\cdots v^{l_n}_n(z):~,} where $v^d_j(z)$ is given by \hwv. For example, in the case $n=2$ and $l=(2,-3)\in\Z^2$, we have $$\omega_l(z) = ~:v^2_1(z) v^{-3}_2(z):~= ~ :\gamma^{x_1}(z)\gamma^{x_1}(z) \beta^{x_2}(z) \beta^{x_2}(z) \beta^{x_2}(z):~.$$ For any conformal structure $L^{\alpha}$ on $\cS(V)$, $\omega_l(z)$ corresponds under the Zhu map to the element $\omega_l\in\cD(V)$ given by \lmclassical. 

\lemma{For each $l\in\Z^n$, the $\cB$-module $\cM_l$ generated by $\omega_l(z)$ is irreducible. Moreover, as a module over $\cB$, \eqn\freemoduleq{\cS(V)= \bigoplus_{l\in\Z^n} \cM_l.}}\thmlab\freemodule
 \proof This is immediate from the description of $\cS^d$ as the irreducible $\cS^0$-module generated by $v_d(z)$, and the fact that $\cS(V)' = \cS'_1\otimes \cdots\otimes \cS'_n$. $\Box$

Note that $\theta^j(z)\circ_0 \omega_l(z) = -l_j \omega_l(z)$, so the $\Z^n$-grading on $\cS(V)$ above is just the eigenspace decomposition of $\cS(V)$ under the family of diagonalizable operators $-\theta^j(z)\circ_0$.

For the remainder of this section, $\gg$ will denote the abelian Lie algebra $$\C^m =  gl(1)\oplus\cdots\oplus gl(1),$$ and $\rho:\gg\ra End(V)$ will be a faithful, diagonal action. Let $A(\rho)\subset \C^n$ be the subspace spanned by $\{\rho(\xi)|~\xi\in\gg\}$. As in the classical setting, we denote $\cS(V)^{\Theta_+}$ by $\cS(V)^{\Theta_+}_{\rho}$ when we need to emphasize the dependence on $\rho$. Clearly $\cS(V)^{\Theta_+}_{\rho} = \cS(V)^{\Theta_+}_{g\cdot \rho}$ for all $g\in GL(m)$, so the family of algebras $\cS(V)^{\Theta_+}_{\rho}$ is parametrized by the points $A(\rho)\in Gr(m,n)$.

Choose a basis $\{\xi^1,\dots,\xi^m\}$ for $\gg$ such that the corresponding vectors $$\rho(\xi^i) = a^i = (a^i_1,\dots,a^i_n)\in\C^n$$ form an orthonormal basis for $A = A(\rho)$.  Let $\theta^{\xi_i}(z)$ be the vertex operator corresponding to $\rho(\xi^i)$, and let $\Theta$ be the subalgebra of $
\cB$ generated by $\{\theta^{\xi_i}(z)|~i=1,\dots,m\}$. By \deftheta, we have 
$$\theta^{\xi_i}(z)=  \sum_{j=1}^n a_j \theta^j(z) = -\sum_{j=1}^n a_j :  \gamma^{x'_j}(z)\beta^{x_j} (z):~.$$ Clearly $\theta^{\xi_i}(z) \theta^{\xi_j}(w)\sim -\bra a^i,a^j\ket (z-w)^{-2} = \delta_{i,j} (z-w)^{-2}$.

If $m<n$, extend the set $\{a^1,\dots,a^m\}$ to an orthonormal basis for $\C^n$ by adjoining vectors $b^i = (b^i_1,\dots,b^i_n)\in \C^n$, for $i=m+1,\dots,n$. Let $$\phi^i(z) = \sum_{j=1}^n b^i_j \theta^j(z) = -\sum_{j=1}^n b^i_j :\gamma^{x'_j}(z)\beta^{x_j}(z):$$ be the corresponding vertex operators, and let $\Phi$ be the subalgebra of $\cB$ generated by $\{\phi^i(z)|~i=m+1,\dots,n\}$. The OPEs
$$\phi^i(z) \phi^j(w)\sim -\bra b^i,b^j\ket (z-w)^{-2},~~~~~~~\theta^{\xi_i}(z)\phi^j(w)\sim -\bra a^i,b^j\ket (z-w)^{-2}$$ show that the $\phi^i(z)$ pairwise commute and each generates a Heisenberg algebra of central charge $-1$, and that $\Phi\subset \cS(V)^{\Theta_+}$. In particular, we have the decomposition $$\cH^1\otimes \cdots \otimes \cH^n = \Theta\otimes \Phi.$$
Next, let $\cW$ denote the subalgebra of $\cB$ generated by $\{L^j(z),W^j(z)|~j=1,\dots,n\}$. Theorem \main~shows that $\cW$ commutes with both $\Theta$ and $\Phi$, so we have the decomposition \eqn\decompofb{\cB = \cW\otimes \Theta\otimes \Phi.} In particular, the subalgebra $\cB' = \cW\otimes \Phi$ lies in the commutant $\cS(V)^{\Theta_+}$. Let $\cM'_l$ denote the $\cB'$-submodule of $\cM_l$ generated by $\omega_l(z)$, which is clearly irreducible as a $\cB'$-module.

In order to describe $\cS(V)^{\Theta_+}$, we first describe the larger space $\cS(V)^{\Theta_>}$ which is annihilated by $\theta^{\xi_i}(k)$ for $i=1,
\dots,m$ and $k>0$. Then $\cS(V)^{\Theta_+}$ is just the subspace of $\cS(V)^{\Theta_>}$ which is annihilated by $\theta^{\xi_i}(0)$, for $i=1,\dots,m$. It is clear from \decompofb~and the irreducibility of $\cM_l$ as a $\cB$-module that $\cS(V)^{\Theta_>}\cap \cM_l = \cM'_l$, so 
 \eqn\msvai{\cS(V)^{\Theta_>} = \bigoplus_{l\in\Z^n} \cM'_l.}

\theorem{As a module over $\cB'$, \eqn\msvaii{\cS(V)^{\Theta_+} = \bigoplus_{l\in A^{\perp}\cap\Z^n} \cM'_l.}}
\proof Let $\omega(z)\in\cS(V)^{\Theta_+}$. Since $\omega$ lies in the larger space $\cS(V)^{\Theta_>}$ which is a direct sum of irreducible, cyclic $\cB'$-modules $\cM'_l$ with generators $\omega_l(z)$, we may assume without loss of generality that $\omega(z) = \omega_l(z)$ for some $l$. An OPE calculation shows that \eqn\latticeope{\theta^{\xi_i}(z)\omega_l(w)\sim -\bra a^i,l \ket \omega_l(w) (z-w)^{-1}.} Hence $\omega_l \in \cS(V)^{\Theta_+}$ if and only if $l$ lies in the sublattice $A^{\perp}\cap \Z^n$. $\Box$

Our next step is to find a {\it finite} generating set for $\cS(V)^{\Theta_+}$. Generically, $A^{\perp}\cap \Z^n$ has rank zero, so $\cS(V)^{\Theta_+} = \cB'$, which is (strongly) generated by the set $$\{\phi^i(z), L^j(z), W^j(z)| i=m+1,\dots,n,~j=1,\dots,n\}.$$ If $A^{\perp}\cap \Z^n$ has rank $r$ for some $0<r\leq n-m$, choose a basis $\{l^1,\dots,l^r\}$ for $A^{\perp}\cap \Z^n$. We claim that for any $l\in A^{\perp}\cap \Z^n$, $\omega_l(z)$ lies in the vertex subalgebra generated by $$\{\omega_{l^1}(z),\cdots,\omega_{l^r}(z),\omega_{-l^1}(z),\dots,\omega_{-l^r}(z)\}.$$ It suffices to prove that given lattice points $l = (l_1,\dots,l_n)$ and $l' = (l'_1,\dots,l'_n)$ in $\Z^n$, $\omega_{l+l'}(z) = k \omega_{l}(z)\circ_d \omega_{l'}(z)$ for some $k\neq 0$ and $d\in\Z$.

First, consider the special case where $l = (l_1,0,\dots,0)$ and $l' = (l'_1,0,\dots,0)$. If $l_1l'_1\geq 0$, we have $\omega_l(z)\circ_{-1}\omega_{l'}(z) = \omega_{l+l'}(z)$. Suppose next that $l_1<0$ and $l'_1>0$, so that $\omega_l(z) =  \beta^{x_1}(z)^{-l_1}$ and $\omega_{l'}(z) = \gamma^{x'_1}(z)^{l'_1}$. Let $$d_1 = min\{ -l_1,l'_1\},~~~~~~e_1 = max\{-l_1,l'_1\},~~~~~~ d = d_1 -1.$$ An OPE calculation shows that 
\eqn\latticecontractioni{\omega_l(z)\circ_d \omega_{l'}(z) = \frac{e_1!}{(e_1-d_1)!} \omega_{l+l'}(z),} where as usual $0! = 1$. Similarly, if $l_1>0$ and $l'_1<0$, we take $d_1 = min\{ l_1,-l'_1\}$,  $e_1 = max\{l_1,-l'_1\}$, and $d = d_1-1$. We have
\eqn\latticecontractionii{\omega_l(z)\circ_d \omega_{l'}(z) = -\frac{e_1!}{(e_1-d_1)!} \omega_{l+l'}(z).}
Now consider the general case $l = (l_1,\dots,l_n)$ and $l' = (l'_1,\dots,l'_n)$. 
For $j=1,\dots,n$, define 
$$d_j = \bigg\{ \matrix{0 & l_j l'_j\geq 0 \cr  min\{ | l_j |, |l'_j|\}, &  l_j l'_j <0 }~~,~~~~~~ e_j = \bigg\{ \matrix{0 & l_j l'_j\geq 0 \cr  max\{ | l_j |, |l'_j|\}, &  l_j l'_j <0 }~,$$
$$k_j= \bigg\{ \matrix{0 & l_j\leq 0\cr  d_j &  l_j >0 }~~,~~~~~~~~d = -1 +\sum_{j=1}^n d_j~~.$$ 
Using \latticecontractioni~and \latticecontractionii~repeatedly, we calculate
$$\omega_l(z)\circ_d \omega_{l'}(z) = \bigg(\prod_{j=1}^n (-1)^{k_j}\frac{e_j!}{(e_j-d_j)!} \bigg)\omega_{l+l'}(z),$$ which shows that $\omega_{l+l'}(z)$ lies in the vertex algebra generated by $\omega_l(z)$ and $\omega_{l'}(z)$. Thus we have proved

\theorem{Let $\{l^1,\dots,l^r\}$ be a basis for the lattice $A^{\perp}\cap \Z^n$, as above. Then $\cS(V)^{\Theta_+}$ is generated as a vertex algebra by $\cB'$ together with the additional vertex operators
$$\omega_{l^1}(z),\dots, \omega_{l^r}(z),~~~\omega_{-l^1}(z),\dots, \omega_{-l^r}(z).$$ In particular, $\cS(V)^{\Theta_+}$ is finitely generated as a vertex algebra.}\thmlab\finitegeneration

In the generic case where $A^{\perp}\cap \Z^n=0$ and $\cS(V)^{\Theta_+} = \cB'$, we claim that $\cS(V)^{\Theta_+}$ has a natural $(n-m)$-parameter family of conformal structures for which the generators $\phi^i(z), L^j(z),W^j(z)$ are primary of conformal weights $1,2,3$, respectively. Note first that $\cW$ has the conformal structure $L_{\cW}(z) = \sum_{j=1}^n L^j(z)$ of central charge $-2n$.

It is well known that for $k\neq 0$ and $c\in\C$, the Heisenberg algebra $\cH$ of central charge $k$ admits a Virasoro element $L^{c}(z) = \frac{1}{2k} j(z)j(z) + c \partial j(z)$ of central charge $1-12 c ^2 k$, under which the generator $j(z)$ is primary of weight one. Hence given $\lambda = (\lambda_{m+1},\dots,\lambda_{n})\in\C^{n-m}$ the Heisenberg algebra generated by $\phi^i(z)$ has a conformal structure $$L^{\lambda_i}(z) = -\frac{1}{2} :\phi^i(z)\phi^i(z):~ + ~\lambda_i \partial \phi^i(z)$$ of central charge $1+12\lambda_i^2$. Since $\phi^i(z)$ and $\phi^j(z)$ commute for $i\neq j$, it follows that $L^{\lambda}_{\Phi}(z) = \sum_{i=m+1}^{n} L^{ \lambda_i}(z)$ is a conformal structure on $\Phi$ of central charge $\sum_{i=m+1}^{n} 1+12 \lambda_i^2$. Finally, $$L_{\cB'}(z) = L_{\cW}(z)\otimes 1 + 1\otimes L^{\lambda}_{\Phi}(z)\in\cW\otimes \Phi = \cB'$$ is a conformal structure on $\cB'$ of central charge $-2n + \sum_{i=m+1}^{n} 1+12\lambda_i^2$ with the desired properties.

When the lattice $A^{\perp}\cap \Z^n$ has positive rank, the vertex algebras $\cS(V)^{\Theta_+}$ have a very rich structure which depends sensitively on $A^{\perp}\cap\Z^n$. In general, the set of generators for $\cS(V)^{\Theta_+}$ given by Theorem \finitegeneration~will not be a set of {\it strong} generators, and the conformal structure $L_{\cB'}$ on $\cB'$ will not extend to a conformal structure on all of $\cS(V)^{\Theta_+}$. 

\theorem{For any action of $\gg$ on $V$, $Com(\cS(V)^{\Theta_+},\cS(V)) = \Theta$. Hence $\cS(V)^{\Theta_+}$ and $\Theta$ form a Howe pair inside $\cS(V)$.}\thmlab\hpva

\proof Since $\cB'\subset \cS(V)^{\Theta_+}$, we have $\Theta\subset Com(\cS(V)^{\Theta_+},\cS(V)) \subset Com(\cB',\cS(V))$, so it suffices to show that $Com(\cB',\cS(V)) = \Theta$. Recall that $\cB' = \cW\otimes \Phi$ and $\Theta\otimes \Phi = \cH^1\otimes\cdots\otimes \cH^n$. Since $Com(\cW^i,\cS_i) = \cH^i$ by Theorem \main, it follows that $Com(\cW,\cS(V)) = \Theta \otimes \Phi$. Then $$Com(\cB',\cS(V)) = Com \big(\Phi, Com(\cW,\cS(V))\big)=Com(\Phi,\Theta\otimes \Phi) = \Theta\otimes Com(\Phi,\Phi) = \Theta.~~ \Box$$

This result shows that we can always recover the action of $\gg$ (up to $GL(m)$-equivalence) from $\cS(V)^{\Theta_+}$, by taking its commutant inside $\cS(V)$. This stands in contrast to Theorem \howepairclassical, which shows that we can reconstruct the action from $\cD(V)^{\gg}$ only when $A^{\perp}\cap \Z^n$ has rank $n-m$. 

\theorem{For any action of $\gg$ on $V$, $\cS(V)^{\Theta_+}$ is a simple vertex algebra.}\thmlab\simplicity
\proof Given a non-zero ideal $\cI\subset\cS(V)^{\Theta_+}$, we need to show that $1\in\cI$. Let $\omega(z)$ be a non-zero element of $\cI$. Since each $\cM'_l$ is irreducible as a module over $\cB'$, we may assume without loss of generality that \eqn\simpleva{\omega(z) = \sum_{l\in\Z^n} c_l \omega_l(z)} for constants $c_l\in\C$, such that $c_l \neq 0$ for only finitely many values of $l$.

For each lattice point $l = (l_1,\dots,l_n)\in\Z^n$, both $\omega_l(z)$ and $\omega_{-l}(z)$ have degree $d = \sum_{j=1}^n |l_j|$ as polynomials in the variables $\beta^{x_j}(z)$ and $\gamma^{x'_j}(z)$. Let $d$ be the maximal degree of terms $\omega_l(z)$ appearing in \simpleva~with non-zero coefficient $c_l$, and let $l$ be such a lattice point for which $\omega_l(z)$ has degree $d$. An OPE calculation shows that

\eqn\calcsimple{\omega_{-l}(z)\circ_{d-1}\omega_{l'}(z)= \bigg\{ \matrix{0 & l'\neq l \cr \bigg(\prod_{j=1}^n (-1)^{k_j} |l_j| ! \bigg)1 & l'=l}} where $k_j =  min\{0,l_j\}$, for all lattice points $l'$ appearing in \simpleva~with non-zero coefficient. It follows from \calcsimple~that
$$\frac{1}{c_l \bigg(\prod_{j=1}^n (-1)^{k_j} |l_j| ! \bigg)} \omega_{-l}(z)\circ_{d-1}\omega(z) = 1.~\Box$$

\subsec{The map $\pi:\cS(V)^{\Theta_+}\ra \cD(V)^{\gg}$}
Equip $\cS(V)$ with the conformal structure $L^{\alpha}$ given by \virbetagamma, for some $\alpha = (\alpha_1,\dots,\alpha_n)\in\C^n$. Suppose first that $A^{\perp}\cap \Z^n$ has rank zero, so that $\cS(V)^{\Theta_+} = \cB'$, and $\cD(V)^{\gg} = \C[e_1,\dots,e_n] = E$. Let $\pi:\cS(V)^{\Theta_+}\ra \cD(V)^{\gg}$ be the map given by \commdiagva. By Lemma \zhugeni, for $j=1,\dots,n$ we have
$$\pi(L^{j}(z)) =  \frac{1}{2}(e_j^2 + e_j),~~~~~~~ \pi(W^j(z)) = \frac{2}{3\sqrt{6}} e_j^3 + \frac{1}{\sqrt{6}} e_j^2 + \frac{1}{3\sqrt{6}} e_j.$$ 
Moreover, \pizhubg~shows that $\pi(\phi^{i}(z)) =  \bra b^i,\alpha\ket - \sum_{j=1}^n b^i_j (e_j +1 )$. Since $\cB'$ is strongly generated by $\{\phi^i(z),L^j(z),W^j(z)|~i=m+1,\dots,n,~j=1,\dots,n\}$, it follows from Lemma \zhugeneration~that $Im(\pi)$ is generated by the collection $$\{\pi(\phi^i(z)), \pi(L^j(z)),\pi(W^j(z))|~i=m+1,\dots,n,~~j=1,\dots,n\}.$$
The map $\pi$ is not surjective, but $Coker (\pi)$ is generated as a module over $Im(\pi)$ by the collection $\{t^{\xi_i}|~i=1,\dots,m\}$, where $t^{\xi_i}$ is the image of $$\pi_{Zh}(\theta^{\xi_i}(z)) = \bra a^i,\alpha\ket - \sum_{j=1}^n a^i_j (e_j+1)$$ in $Coker(\pi) = E/\pi(\cB')$. Unlike the case where $V$ is one-dimensional, $\pi$ depends on the choice of $\alpha$.

Suppose next that the lattice $A^{\perp}\cap \Z^n = 0$ has positive rank. Clearly $\pi_{Zh}(\cM_l) = M_l$ for all $l$, so $\pi(\cM'_l)\subset M_l$. This map need not be surjective, but since $M_l$ is the free $E$-module generated by $\omega_l$, and $E/ \pi(\cB')$ is generated as a $\pi(\cB')$-module by $\{t^{\xi_i}|~i=1,\dots,m\}$, it follows that each $M_l/\pi(\cM'_l)$ is generated as a $\pi(\cB')$-module by $\{ t^{\xi_i}_l|~i=1,\dots,m\}$, where $t^{\xi_i}_l$ is the image of $\pi_{Zh}(\theta^{\xi_i}(z))\omega_l$ in $M_l/\pi(\cM'_l)$.

\theorem{For any action of $\gg$ on $V$, $Coker (\pi)$ is generated as a module over $Im(\pi)$ by the collection $\{t^{\xi_i}|~i=1,\dots,m\}$. In particular, $Coker(\pi)$ is a finitely generated module over $Im(\pi)$ with generators corresponding to central elements of $\cD(V)^{\gg}$.}

\proof First, since $\pi(\omega_l(z)) = \omega_l$ for all $l$, it is clear that the generators $t^{\xi_i}_l$ of $M_l / \pi(\cM'_l)$ lie in the $Im(\pi)$-module generated by $\{t^{\xi_i}|~i=1,\dots,m\}$, which proves the first statement. Finally, the fact that the elements $\pi_{Zh}(\theta^{\xi_i}(z))$ corresponding to $t^{\xi_i}$ each lie in the center of $\cD(V)^{\gg}$ is immediate from \perpcond. $\Box$

\subsec{A vertex algebra bundle over the Grassmannian $Gr(m,n)$}
As $\rho$ varies over the space $R^0(V)$ of effective actions, recall that $\cS(V)^{\Theta_+}_{\rho}$ is uniquely determined by the point $A(\rho)\in  Gr(m,n)$. The algebras $\cS(V)^{\Theta_+}_{\rho}$ do not form a fiber bundle over $Gr(m,n)$. However, the subspace of $\cS(V)^{\Theta_+}_{\rho}$ of degree zero in the $A(\rho)^{\perp}\cap \Z^n$-grading \msvaii~is just $\cB'_{\rho} = \cB'$, and the algebras $\cB'_{\rho}$ form a bundle of vertex algebras $\cE$ over $Gr(m,n)$. The classical analogue of $\cE$ is not interesting; it is just the trivial bundle whose fiber over each point is the polynomial algebra $E$. 

For each $\rho$, recall that $\cB'_{\rho} = \cW_{\rho}\otimes \Phi_{\rho}$, where $\cW_{\rho}$ is generated by $\{L^j(z),W^j(z)|~j=1,\dots,n\}$, and $\Phi_{\rho}$ is generated by $\{\phi^i(z)|~i=m+1,\dots,n\}$. Since $\cW_{\rho}$ is independent of $\rho$, it gives rise to a trivial subbundle of $\cE$. As a vector space, note that $\Phi_{\rho} = Sym \big(\bigoplus_{k\geq 1} A(\rho)^{\perp}_k\big)$, where $A(\rho)^{\perp}_k$ is the copy of $A(\rho)^{\perp}$ spanned by the vectors $\partial^k \phi^i(z)$ for $i=m+1,\dots,n$. It follows that the factor $\Phi_{\rho}$ in the fiber over $A(\rho)$ gives rise to the following subbundle of $\cE$: \eqn\vabundle{Sym\big(\bigoplus_{k\geq 1} \cF_k\big),} where $\cF_k$ is the quotient of the rank $n$ trivial bundle over $Gr(m,n)$ by the tautological bundle. Since each $\cF_k$ has weight $k$, the weighted components of the bundle \vabundle~are all finite-dimensional. The non-triviality of this bundle is closely related to Theorem \hpva.

\newsec{Vertex algebra operations and transvectants on $\cD(V)^{\gg}$}

If we fix a basis $\{x_1,\dots,x_n\}$ for $V$ and a dual basis $\{x'_1,\dots, x'_n\}$ for $V^*$, $\cS(V)$ has a basis consisting of iterated Wick products of the form $$\mu(z) = ~:\partial^{k_1} \gamma^{x'_{i_1}}(z)\cdots \partial^{k_r}\gamma^{x'_{i_r}}(z)\partial^{l_1}\beta^{x_{j_1}}(z)\cdots \partial^{l_s}\beta^{x_{j_s}}(z):~.$$ Define gradings {\it degree} and {\it level} on $\cS(V)$ as follows: $$ deg (\mu) = r+s,~~~~~~ lev(\mu) =\sum_{i=1}^r k_i + \sum_{j=1}^s l_j,$$ and let $\cS(V)^{(n)}[d]$ denote the subspace of level $n$ and degree $d$. The gradings \eqn\grd{\cS(V) = \bigoplus_{n\geq 0} \cS(V)^{(n)} = \bigoplus_{n,d \geq 0} \cS(V)^{(n)}[d] = \bigoplus_{d \geq 0} \cS(V)[d] } are clearly independent of our choice of basis on $V$, since an automorphism of $V$ has the effect of replacing $\beta^{x_i}$ and $\gamma^{x'_i}$ with linear combinations of the $\beta^{x_i}$'s and $\gamma^{x'_i}$'s, respectively. 

Let $\sigma:\cD(V)\ra gr\cD(V) = Sym(V\oplus V^*)$ be the map \eqn\mapsigma{x'_{i_1}\cdots x'_{i_r}\frac{\partial}{\partial x'_{j_1}}\cdots\frac{\partial}{\partial x'_{j_s}}\mapsto  x'_{i_1}\cdots x'_{i_r} x_{j_1}\cdots x_{j_s},} which is a linear isomorphism. Any bilinear product $*$ on $Sym(V\oplus V^*)$ corresponds to a bilinear product on $\cD(V)$, which we also denote by $*$, as follows:
$$\omega*\nu = \sigma^{-1}(\sigma(\omega)* \sigma(\omega)),$$ for $\omega,\nu \in \cD(V)$,  Moreover, $\omega_1,\dots,\omega_k$ generate $\cD(V)$ as a ring if and only if $\sigma(\omega_1),\dots,\sigma(\omega_k)$ generate $Sym(V\oplus V^*)$ as a ring. The map $f: Sym(V\oplus V^*)\ra \cS(V)^{(0)}$ given by \eqn\mapf{x'_{i_1}\cdots x'_{i_r} x_{j_1}\cdots x_{j_s}, \mapsto  ~: \gamma^{x'_{i_1}}(z)\cdots \gamma^{x'_{i_r}}(z)\beta^{x_{j_1}}(z)\cdots\beta^{x_{j_s}}(z):~,} is a linear isomorphism, so that $f\circ\sigma: \cD(V)\ra \cS(V)^{(0)}$ is a linear isomorphism as well.

$\cS(V)^{(0)}$ has a family of bilinear products $*_k$ which are induced by the circle products on $\cS(V)$. Given $\omega(z),\nu(z)\in \cS(V)^{(0)}$, define \eqn\starn{\omega(z)*_k \nu(z) = p (\omega(z)\circ_k \nu(z)),} where $p:\cS(V)\ra \cS(V)^{(0)}$ is the projection onto the subspace of level zero. 
Clearly $\omega(z) *_k \nu(z) = 0$ whenever $k<-1$ because $p\circ \partial$ acts by zero on $\cS(V)^{(0)}$. For $k\geq -1$, $*_{k}$ is homogeneous of degree $-2k-2$.

Via \mapf, we may pull back the products $*_k$, $k\geq -1$ to obtain a family of bilinear products on $Sym(V\oplus V^*)$, which we also denote by $*_k$. In fact, these products have a classical description. Let \eqn\defgamma{\Gamma = \sum_{i=1}^n \frac{\partial}{\partial x_i}\otimes \frac{\partial}{\partial x'_i} - \frac{\partial}{\partial x'_i}\otimes \frac{\partial}{\partial x_i},} and define the $k$th transvectant$^1$ \footnote{}{I thank N. Wallach for explaining this construction to me.}on $Sym(V\oplus V^*)$ by $$[,]_k: Sym(V\oplus V^*)\otimes Sym(V\oplus V^*)\ra Sym(V\oplus V^*),~~~~ [\omega,\nu]_k = m\circ \Gamma^k (\omega\otimes \nu).$$ Here $m$ is the multiplication map sending $\omega\otimes \nu\mapsto \omega\nu$. 

\theorem{The product $*_k$ on $Sym(V\oplus V^*)$ given by \starn~coincides with the transvectant $[,]_{k+1}$ for $k\geq -1$.}\thmlab\transvectant

\proof First consider the case $k=-1$. In this case $[,]_0$ is just ordinary multiplication. Recall the formula  $$:(:ab:)c:-:abc:=\sum_{k\geq0}{1\over(k+1)!}\left(:(\partial^{k+1}a)(b\circ_k c): +(-1)^{|a||b|}:(\partial^{k+1}b)(a\circ_k c):\right),$$ which holds for any vertex operators $a,b,c$ in a vertex algebra $\cA$. It follows that the associator ideal in $\cS(V)$ under the Wick product is annihilated by the projection $p$. Similarly, the commutator ideal in $\cS(V)$ under the Wick product is annihilated by $p$, so $\cS(V)^{(0)}$ is a polynomial algebra with product $*_{-1}$, and $f: Sym(V\oplus V^*)\ra \cS(V)^{(0)}$ is an isomorphism of polynomial algebras. Hence given $\omega,\nu\in Sym(V\oplus V^*)$, we have $[\omega,\nu]_0 = \omega\nu = \omega *_{-1} \nu$.

Next, if $k\geq 0$, it is clear from the definition of the vertex algebra products $\circ_k$ that given $\omega(z),\nu(z)\in\cS(V)^{(0)}$, $\omega(z)*_k\nu(z)$ is just the sum of all possible contractions of $k+1$ factors of the form $\beta^{x_i}(z)$ or $\gamma^{x'_i}(z)$ appearing in $\omega(z)$ with $k+1$ factors of the form $\beta^{x_i}(z)$ or $\gamma^{x'_i}(z)$ appearing in $\nu(z)$. Here the contraction of $\beta^{x_i}(z)$ with $\gamma^{x_j}(z)$ is $\delta_{i,j}$, and the contraction of $\gamma^{x_i}(z)$ with $\beta^{x_j}(z)$ is $-\delta_{i,j}$. Similarly, it follows from \defgamma~that given $\omega,\nu \in Sym(V\oplus V^*)$, $[\omega,\nu]_{k+1}$ is the sum of all possible contractions of $k+1$ factors of the form $x_i$ or $x'_i$ appearing in $\omega$ with $k+1$ factors of the form $x_i$ or $x'_i$ appearing in $\nu$. The contraction of $x_i$ with $x'_j$ is $\delta_{i,j}$ and the contraction of $x'_i$ with $x_j$ is $-\delta_{i,j}$. Since $f: Sym(V\oplus V^*)\ra \cS(V)^{(0)}$ is the algebra isomorphism sending $x_i\mapsto \beta^{x_i}(z)$ and $x'_i\mapsto \gamma^{x'_i}(z)$, the claim follows. $\Box$

Via $\sigma:\cD(V)\ra Sym(V\oplus V^*)$ the products $*_k$ on $Sym(V\oplus V^*)$ pull back to bilinear products on $\cD(V)$, which we also denote by $*_k$. These products satisfy $\omega *_k \nu \in \cD(V)_{(r+s-2k-2)}$ for $\omega\in\cD(V)_{(r)}$ and $s\in \cD(V)_{(s)}$. It is immediate from Theorem \transvectant~that $*_{-1}$ and $*_0$ correspond to the ordinary associative product and bracket on $\cD(V)$, respectively. Since the circle product $\circ_0$ is a derivation of every $\circ_k$, it follows that $\omega*_0$ is a derivation of $*_k$ for all $\omega\in \cD(V)$ and $k\geq -1$.

We call $\cD(V)$ equipped with the products $\{*_k|~k\geq -1\}$ a $*$-algebra. A similar construction goes through in other settings as well. For example, given a Lie algebra $\gg$ equipped with a symmetric, invariant bilinear form $B$, $\gU\gg$ has a $*$-algebra structure (which depends on $B$). Given a $*$-algebra $\cA$, we can define $*$-subalgebras, $*$-ideals, quotients, and homomorphisms in the obvious way. If $V$ is a module over a Lie algebra $\gg$, $\cD(V)^{\gg}$ is a $*$-subalgebra of $\cD(V)$ because the action of $\xi\in\gg$ is given by $[\tau(\xi),-] = \tau(\xi)*_0$ which is a derivation of all the other products. 

Given elements $\omega_1,\dots,\omega_k\in \cD(V)^{\gg}$, examples are known where $\omega_1,\dots,\omega_k$  do not generate $\cD(V)^{\gg}$ as a ring, but do generate $\cD(V)^{\gg}$ as a $*$-algebra.$^2$\footnote{}{I thank N. Wallach for pointing this out to me.} This phenomenon occurs in our main example, in which $\gg$ is the abelian Lie algebra $\C^m$ acting diagonally on $V= \C^n$. Recall that $\cD(V)^{\gg} = \bigoplus_{l\in A^{\perp}\cap \Z^n} M_l$, where $M_l$ is the free $E$-module generated by $\omega_l$. Suppose that $A^{\perp}\cap \Z^n$ has rank $r$, and let $\{l^i = (l^i_1,\dots,l^i_n)|~i=1,\dots,r\}$ be a basis for $A^{\perp}\cap\Z^n$. In general, the collection \eqn\notabasis{e_1,\dots,e_n, ~~~\omega_{l^1},\dots,\omega_{l^r},~~~ \omega_{-l^1},\dots,\omega_{-l^r}} is too small to generate $\cD(V)^{\gg}$ as a ring. 

\theorem{$\cD(V)^{\gg}$ is generated as a $*$-algebra by the collection \notabasis. Moreover, $\cD(V)^{\gg}$ is simple as a $*$-algebra.}

\proof To prove the first statement, it suffices to show that given lattice points $l = (l_1,\dots,l_n)$ and $l' = (l'_1,\dots,l'_n)$, $\omega_{l+l'}$ lies in the $*$-algebra generated by $\omega_l$ and $\omega_{l'}$. For $j=1,\dots,n$, define 
$$d_j = \bigg\{ \matrix{0 & l_j l'_j\geq 0 \cr  min\{ | l_j |, |l'_j|\}, &  l_j l'_j <0 }~~,~~~~~~ e_j = \bigg\{ \matrix{0 & l_j l'_j\geq 0 \cr  max\{ | l_j |, |l'_j|\}, &  l_j l'_j <0 }~,$$
$$k_j= \bigg\{ \matrix{0 & l_j\leq 0\cr  d_j &  l_j >0 }~~,~~~~~~~~d = -1 +\sum_{j=1}^n d_j~~.$$ 
The same calculation as in the proof of Theorem \finitegeneration~shows that
$$\omega_l *_d \omega_{l'} = \bigg(\prod_{j=1}^n (-1)^{k_j}\frac{e_j!}{(e_j-d_j)!}\bigg) \omega_{l+l'},$$ which shows that $\omega_{l+l'}$ lies in the $*$-algebra generated by $\omega_l$ and $\omega_{l'}$.

As for the second statement, the argument is analogous to the proof of Theorem \simplicity. Given a non-zero $*$-ideal $I\subset\cD(V)^{\gg}$, we need to show that $1\in I$. Let $\omega$ be a non-zero element of $I$. It is easy to check that for $i,j=1,\dots, n$, and $l\in A^{\perp}\cap \Z^n$, we have
$$e_i *_1 e_j = - \delta_{i,j},~~~~~e_i *_1 \omega_l = 0$$
By applying the operators $e_i *_1$ for $i=1,\dots,n$, we can reduce $\omega$ to the form 
\eqn\simpleclassical{\sum_{l\in\Z^n} c_l \omega_l} for constants $c_l\in\C$, such that $c_l \neq 0$ for only finitely many values of $l$. We may assume without loss of generality that $\omega$ is already of this form. Let $d$ be the maximal degree (in the Bernstein filtration) of terms $\omega_l$ appearing in \simpleclassical~with non-zero coefficient $c_l$, and  let $l$ be such a lattice point for which $\omega_l$ has degree $d$. We have
$$\omega_{-l} *_{d-1}\omega_{l'}= \bigg\{ \matrix{0 & l'\neq l \cr \bigg(\prod_{j=1}^n (-1)^{k_j}|l_j| ! \bigg)1 & l'=l}$$ where $k_j = min\{0,l_j\}$, for all $l'$ appearing in \simpleclassical. Hence
$$\frac{1}{c_l \bigg(\prod_{j=1}^n (-1)^{k_j} |l_j| ! \bigg)} \omega_{-l} *_{d-1}\omega = 1.~\Box$$

\footatend\vfill\supereject\immediate\closeout\rfile\writestoppt
\baselineskip=14pt\centerline{{\bf References}}\bigskip{\frenchspacing%
\parindent=20pt\escapechar=` \input refs.tmp\vfill\eject}\nonfrenchspacing

\end